\documentclass[11pt]{article}

\usepackage{latexsym}
\usepackage{amsmath}
\usepackage{amsfonts}
\usepackage{amssymb}
\usepackage{graphpap}
\usepackage{epsfig}
\usepackage{amscd}
\usepackage{amsthm}
\usepackage{enumerate}
\usepackage{eucal} 

\font\got=eufm10 at 11pt

\newcommand{\C}[0]{{\mathbb C}}

\newcommand{\e}[0]{\varepsilon}
\newcommand{\f}[0]{\varphi}

\newcommand{\LL}[0]{{\cal L}}

\newcommand{\N}[0]{{\mathbb N}}
\newcommand{\h}[0]{\text{\got h}}
\newcommand{\He}[0]{{\mathbb H}}
\newcommand{\pd}[0]{\partial}

\newcommand{\R}[0]{\mathbb{R}}

\newcommand{\z}[0]{\mathcal A}

\newcommand{\cH}[0]{{\mathcal H}}
\newcommand{\cK}[0]{{\mathcal K}}

\newcommand{\nnor}[1]{\|\hskip -4.4pt\| #1 \|\hskip -4.4pt\|}
\newcommand{\nor}[1]{\| #1 \|}
\newcommand{\Nor}[1]{\left |\hskip -0.7pt\left | #1 \right |\hskip -0.7pt\right |}

\newcommand{\mn}[2]{\{ #1\, ;#2 \}}
\newcommand{\sk}[2]{\langle #1 , #2\rangle}
\newcommand{\Sk}[2]{\Big\langle #1 , #2\Big\rangle}

\newcommand{\rim}[1]{\underline{\overline{\rm #1}}}

\newcommand{\dok}{\noindent {\bf Proof.}\ }

\newtheorem{thm}{Theorem}
\newtheorem{lema}{Lemma}
\newtheorem{prop}{Proposition}
\newtheorem{cor}{Corollary}

\title{Linear dimension-free estimates for the Hermite-Riesz transforms
\footnote{The first author was 
supported by the Ministry of Higher Education, Science and Technology of Slovenia (research program Analysis and Geometry, contract no. P1-0291). The second author was supported by NSF grant DMS 0501067.}}
\author{
Oliver Dragi\v{c}evi\'c 
and 
Alexander Volberg
}
\date{November 5, 2008}

\begin{document}

\setcounter{page}{1}

\maketitle
\begin{abstract}
We utilize the Bellman function technique to prove a bilinear dimension-free inequality for the Hermite operator. The Bellman technique is applied here to a non-local operator, which at first did not seem to be feasible. 
As a consequence of our bilinear inequality one proves dimension-free boundedness for the Riesz-Hermite transforms on $L^p$ with linear growth in terms of $p$. A feature of the proof is a theorem establishing $L^p(\R^n)$ estimates for a class of spectral multipliers with bounds independent of $n$ and $p$.
Connections with 
known
results on the Heisenberg group 
as well as with results for Hilbert transform along the parabola
are also explored.
We believe 
our approach is quite universal in the sense that one could apply it to a whole range of Riesz transforms arising from various differential operators. 
As a first step towards this goal we prove our dimension-free bilinear embedding theorem for quite a general family of Schr\"odinger semigroups.


\end{abstract}
%
\section{Introduction}
\label{intro}

The purpose of this paper is to 
prove
bilinear $L^p\times L^q\rightarrow L^1$ embeddings associated with the Hermite operator on $\R^n$.
Similar results were proved in \cite{DV} for the usual Laplacian and the Ornstein-Uhlenbeck operator. 
In this paper
we manage to apply the same method (Bellman functions) and treat operators with potentials, such as the Hermite Laplacian, which is a
novelty compared to \cite{DV}. All of our embedding theorems are
dimension-free and exhibit linear estimates in terms of $p$. This feature is due to special properties of
the concrete Bellman function we use. It is further exploited in
\cite{DV1} where the embedding theorem is proved for extensions
generated by second-order Schr\"odinger-type operators in
divergence form with real coefficients (i.e. regarding Kato's problem with real
matrix).

As the main yet simple consequence of our embedding theorem we
give dimension-free $L^p$ estimates of Riesz transforms for the
Hermite semigroup. Our approach already gave such estimates for
euclidean and Ornstein-Uhlenbeck semigroups, see \cite{DV}, and
there we announced analogous results for a wider class of differential
operators, provided only that their spectral properties are not
``too singular", so that we can construct an efficient
dimension-free passage from the embedding theorems. In the Hermite
case the latter appear here, to the best of our knowledge, for the first time in literature.

Regarding the applications to the 
Riesz transforms we should say that for euclidean,
Ornstein-Uhlenbeck 
and Hermite semigroups 
our approach offers a unified way of treating various 
vector-valued operators of this kind 
and obtaining dimension-free estimates.
Apart from this uniformity, it yields explicit estimates of the $L^p$ norms for Riesz transforms in terms of $p$. In the case considered here -- the one of the Hermite operator -- we establish linear behaviour with respect to $p$. This result seems to be new and we believe it to be sharp.
We also prove
dimension-free estimates of iterations of Riesz transforms. 
In order to yield such estimates we present 
bounds on $L^p(\R^n)$ 
for spectral multipliers arising from functions analytic at infinity.
These bounds are absolute, i.e. independent of $n$ and $p$.
Again, we believe that exactly the same treatment can be applied
to a wide class of differential operators (infinitesimal
generators), under the condition that their spectral properties are not
``too singular".

We tried to make the paper as self-contained as possible, given that we were naturally obliged to refer to some basic facts.

\subsection*{Statement of main results}
Take $p\in (1,\infty)$ and denote by $q$ its conjugate exponent $p/(p-1)$.
We will use the notation $p^*:=\max\{p,q\}$.

The {\it Hermite operator} $L$ is, for a test function $u$ on $\R^n$, defined as
$$
Lu(x)=-\Delta u(x)+|x|^2u(x)\,.
$$
An equivalent description of $L$ is given by
\begin{equation}
\label{nikolaeva}
L=\frac12\sum_{j=1}^n(\z_j\z_j^*+\z_j^*\z_j)\,,
\end{equation}
where
\begin{equation}
\label{mole}
\z_j =-\frac{\pd}{\pd x_j} + x_j
\ \hskip 15pt
\text{and}
\ \hskip 15pt
\z_j^* =\frac{\pd}{\pd x_j} + x_j
\end{equation}
are the {\it creation} and {\it annihilation} operator, respectively.
The operator $L$ is positive, meaning that $\sk{Lu}{u}\geqslant0$.
For a thorough discussion of $L$ we refer the reader to \cite{T}.

We will be dealing with the operator semigroup $\{P_t :=e^{-t\sqrt{L}}\}_{t>0}$.
Denote
$
\widetilde u(x,t)  =  P_t u(x)\,.
$
This function solves on $\R^n\times (0,\infty)$ the differential equation
\begin{equation}
\label{ddd_ddd}
\Big(\frac{\pd^2 }{\pd t^2}-L\Big)\widetilde u = 0,
\end{equation}
with
$
\widetilde u(x,0) =u(x)\,.
$

For a given smooth $\C^N-$valued function  $\phi=(\phi_1(x,t),\hdots,\phi_N(x,t))$ on $\R^n\times (0,\infty)$ denote
\begin{equation}
\label{debeluh}
\aligned
\nor{\phi}_*^2
& =  \Big|\frac{\pd \phi}{\pd t}\Big|^2 +
\sum_{j=1}^n \Big|\frac{\pd \phi}{\pd x_j}\Big|^2 + |x|^2 |\phi(x,t)|^2\\
& =\Big|\frac{\pd \phi}{\pd t}\Big|^2  + \frac12 \sum_{j=1}^n (|\z_j \phi(x,t)|^2 + |\z_j^* \phi(x,t)|^2)\,.
\endaligned
\end{equation}
Here, as usual,
$$
\frac{\pd\phi}{\pd x_j}=\left(\frac{\pd\phi_1}{\pd x_j},\hdots,\frac{\pd\phi_N}{\pd x_j}\right)\,.
$$
The $L^p$ norm of a $\C^N-$valued test function $\psi$ on $\R^n$ is of course
$
\big(\int_{\R^n}\nor{\psi(x)}_{\C^N}^{p}\,dx\big)^{1/p}.
$

We are ready to state the key proposition, which can be considered as the bilinear embedding theorem (or bilinear Littlewood-Paley theorem) for the Hermite extension.

\begin{thm}
\label{bilinher}
There is an absolute constant $C>0$ such that for 
arbitrary natural numbers $M,N,n$,  
any pair $f:\R^n\rightarrow\C^M$ and $g:\R^n\rightarrow\C^N$ of $C_c^\infty$ test functions and any
$p>1$ we have
$$
\int_0^\infty\int_{\R^n}\nor{P_t f(x)}_*\nor{ P_t g(x)}_*\,dx\,t\,dt\leqslant C(p^*-1)\nor{f}_p\nor{g}_q\,.
$$
\end{thm}
Note that these estimates are dimension free.

\subsubsection*{Schr\"odinger operators with positive potentials}

The preceding theorem can without changing the statement be generalized to a large class of Schr\"odinger operators $L=-\Delta+V(x)$. Here 
$V$ is a non-negative function on $\R^n$. 
By $P_t$ denote the (Poisson) operator semigroup whose infinitesimal generator is $L^{1/2}$.
Let $\cK_t$ be the heat kernel, i.e. the kernel associated to the semigroup generated by $L$.
Assume the following conditions on $V$:
\begin{enumerate}[(a)]
\item
\label{heatkato}
Kato's inequality
 \begin{equation*}
\cK_t(x,y)\leqslant  Ct^{-\frac{n}{2}}\,{e^{-\frac{a}{t}|x-y|^2}} 
 \end{equation*}
and the kernel being sub-probability, i.e.
 \begin{equation*}
\int_{\R^n} \cK_t(x,y) \,dy\leqslant 1
\,.
 \end{equation*}
for all $x\in\R^n$. Also, $\cK_t$ needs to be non-negative. 

\item
\label{gradest}
Gradient estimates for the heat kernel:
\begin{equation*}
\aligned
\bigg|{\pd \over \pd x_j}\cK_{t}(x,y)\bigg|
& \leqslant C
t^{-{n+1\over 2}}\,e^{-{a\over t}|x-y|^2}
\\
\bigg|\frac{\pd}{\pd t} \cK_t(x,y)\bigg| 
& \leqslant Ct^{-\frac{n}{2}-1}\, e^{-\frac{a}{t}|x-y|^2}\,.
\endaligned
\end{equation*}

\item
\label{lim}
If $g\in C_c^\infty$ 
then
$$
\aligned
& \lim_{t\rightarrow 0}\int_{\R^n} \bigg|t\frac{\pd P_t g}{\pd t}(x)\bigg|\,dx=0\\
& \lim_{t\rightarrow\infty}\int_{\R^n} t\frac{\pd P_t g}{\pd t}(x)\,dx=0\,.
\endaligned
$$

\item
\label{emb}
For any bounded, non-negative, compactly supported function 
$\f$ and some $C_0>0$ which does not depend on $n$,
$$
\int_{\R^n}\int_0^\infty P_t\f(x)\,t\,dt\,V(x)\,dx\leqslant C_0\nor{\f}_1\,.
$$
\end{enumerate}
The conditions which imply \eqref{heatkato} were studied in \cite{Si}, for example. In general there exists a vast literature on estimates of heat kernels.

Under the conditions \eqref{heatkato} -- \eqref{emb} we get exactly the same statement as in Theorem \ref{bilinher}. We emphasize that the constants $C,a$ from the conditions \eqref{heatkato} and \eqref{gradest} are allowed to be arbitrary (meaning also dependent on the dimension $n$), and still they allow us to deduce dimension-free estimates in the analogue of Theorem \ref{bilinher} for general potential $V$ as above.

\subsubsection*{Riesz transforms}

We mentioned above that 
as 
an application
of the embedding theorem 
we obtain 
an $L^p$ estimate of the corresponding {\it Riesz transforms}, introduced in \cite{T} as 
\begin{equation*}
\label{riesz}
R_j = \z_j L^{-\frac12}\,,\hskip 10pt R_j^* = \z_j^* L^{-\frac12}\,,\hskip 10pt j=1,\hdots,n\,.
\end{equation*}
A proof establishing their $L^p$ boundedness can be found in \cite{T}. See also \cite{T3}.
This result has later been enhanced, as we describe now.

Define
$$
{\bf R}
f  =\Big(\sum_{j=1}^n|R_jf|^2+\sum_{j=1}^n|R_j^*f|^2\Big)^{1\over2}\,.
$$
It was 
proved in \cite{HRST} and shortly afterwards in \cite{L-P} that the $L^p$ bounds for ${\bf R}$ admit estimates
from above with constants not depending on the dimension. 
Apparently, however, the authors of \cite{HRST} and \cite{L-P} missed 
the fact that this result can be deduced, by means of transference, from an earlier work \cite{CMZ} addressing analogous questions on the Heisenberg group. 
Moreover,
the approach from \cite{CMZ}
also 
permits explicit information about behaviour of the estimates with respect to $p$; it seems \cite{HRST} and \cite{L-P} do not contain that. Namely,  
the main result of \cite{CMZ} draws a close connection between 
Heisenberg-Riesz transforms and the Hilbert transform along the parabola in $\R^2$. The $L^p$ boundedness of the latter has been known for a long time \cite{SW}, yet only not long ago, owing to several highly nontrivial results due to Seeger, Tao and Wright 
\cite{STW, Tao, TW}, 
has there been major improvement as to the behaviour of estimates in terms of $p$. 
As a result it is possible to obtain, for arbitrary $\e>0$, the 
estimate 
$\nor{{\bf R}f}_p\leqslant C_\e p^{1+\e}\nor{f}_p$.
We shall devote 
section \ref{Heisenberg}
to the purpose of explaining these connections and developments. 

As a consequence of our Theorem \ref{bilinher} we are able to further sharpen this result and obtain linear estimates, which we believe to be optimal.

\begin{cor}
\label{RieszH}
There exists $C>0$ such that 
for 
any $1<p<\infty$, $n\in\N$ and $f\in L^p(\R^n)$, 
\begin{equation}
\label{tuva}
\nor{ 
{\bf R}
f}_p\leqslant 
C(p^*-1)\nor{f}_p\,.
\end{equation}
\end{cor}

This linear inequality either  
gives indirect evidence in favour of the Hilbert transform along parabola admitting linear $L^p$ estimates (see \eqref{podpole} in section \ref{Heisenberg}), or else indicates that the method involving the Bellman function is yet sharper than the passage through the Heisenberg groups described above. We 
believe the second option to be more likely.

\bigskip
The estimate \eqref{tuva} can be generalized to encompass 
iterations of Riesz transforms.
A theorem of the same type was proved recently in \cite{L-P}, but apparently without the numerical estimate of $C(p)$. Let us formulate our result rigorously. 
Write
$$
{\cal C}_d=\{\text{compositions of }d\text{ operators among } R_1,\hdots,R_n,R_1^*,\hdots,R_n^*\} \,.
$$
\begin{cor}
\label{schwertleite}
Assuming conditions and notation as in Corollary \ref{RieszH},
$$
\bigg\|\bigg(\sum_{R^d\in{\cal C}_d}|R^df|^2\bigg)^{1/2}\bigg\|_p\leqslant C^d(p^*-1)^{d}\nor{f}_p\,.
$$
\end{cor}

\subsubsection*{Spectral multipliers}
It is worth mentioning that in order to prove the above corollaries we come up with 
another result which we believe to be of independent interest. It deals with 
$L^p$-boundedness of
spectral multipliers for Hermite expansions. Theorems of this type 
were
obtained in \cite{M} and \cite{T}. 
Our theorem is different in terms of the assumptions 
laid on the multiplier
and also because we obtain absolute bounds for the norms, i.e. such which do not depend on $n$ or $p$. 
Further discussion of the theorem and its proof 
are to be found
in section \ref{sm}.

\begin{thm}
\label{multiplier}
Let function $\Psi$ be analytic at infinity. 
Then for all $p\in [1, \infty]$ and all $n\in\N$, the operator 
$\Psi(L)$
is bounded on $L^p(\R^n)$ by a constant which depends only on $\Psi$, i.e. it does not depend neither on $n$ nor $p$.
\end{thm}



\bigskip
{\bf Acknowledgements.} We are grateful to
Giancarlo Mauceri and Stefano Meda for discussions on 
spectral multipliers for the Hermite operator.
Our deep gratitude goes to Fulvio Ricci, Adam Sikora, Jos\'e Luis Torrea and Fran\c coise Lust-Piquard for many helpful explanations regarding the Hermite semigroup. We also thank Jim Wright for conveying to us recent results about  Hilbert transforms along curves. 

Finally, we wish to thank Centro Di Ricerca Matematica Ennio De Giorgi, Pisa, for the hospitality extended to us for the two weeks in March 2007.

\section{Bellman function}

\label{DSCH}

Throughout the section we  work with $p\geqslant 2$, $q=p/(p-1)$ and $\delta=q(q-1)/8$. Observe that $\delta\sim (p-1)^{-1}$. 

The crucial part in our proofs will be played by the function $Q$, given  by
\begin{equation}
\label{sarmat}
Q(\zeta,\eta,Z,H)=2(Z+H)-|\zeta|^p-|\eta|^q-\delta\tilde
Q(\zeta,\eta)\,,
\end{equation}
where
$$
\widetilde Q(\zeta,\eta)=\left\{
\aligned
& |\zeta|^2|\eta|^{2-q} & ; & \ \ |\zeta|^p\leqslant |\eta|^q\\
& \frac{2}{p}|\zeta|^p+\left(\frac{2}{q}-1\right)|\eta|^q
& ; &\ \ |\zeta|^p\geqslant |\eta|^q\,\ \,.
\endaligned\right.
$$
Note that the definition of $Q$ depends on $p$. 
Such a function, defined on a subdomain in $\R^4$, was
first introduced by F. Nazarov and S. Treil \cite{NT}. Here it is
defined in the domain
$$
\Omega :=\mn{(\zeta,\eta,Z,H)\in\C^M\times\C^N\times\R\times\R}{|\zeta|^p \leqslant Z, |\eta|^q \leqslant H}\,.
$$

Function $Q$ is in $C^1(\Omega)$ and its second derivatives are continuous except on $\gamma:=\{|\zeta|^p=|\eta|^q\}$ but everywhere in $\Omega$ (including $\gamma$) these second derivatives or their one-sided limits can be estimated in a way which suits our purposes very well.

The Hessian matrix of $Q$ is denoted by $d^2Q$. Thus $d^2Q$ is a matrix-valued function
which maps vector $\omega\in\Omega$ into the matrix with entries $\frac{\pd^2Q}{\pd\alpha\pd\beta}(\omega)$,
where $\alpha$ and $\beta$ range over $\zeta_j,\overline{\zeta_j},\eta_k,\overline{\eta_k},Z,H$ for $j=1,\hdots,M$, $k=1,\hdots,N$.

\begin{thm}
\label{bejaz}
Choose $\omega=(\zeta,\eta,Z,H)\in\Omega$. Then
\begin{enumerate}[{\rm (i)}]
\item
\label{anterija}
$Q(\omega)\leqslant 2(Z+H)$.
\end{enumerate}
There exists  $\tau=\tau(|\zeta|,|\eta|)>0$ such that
\begin{enumerate}[{\rm (i)}]
\addtocounter{enumi}{1}
\item
\label{kupres}
$
{-d^2 Q(\omega) \geqslant \delta\big(\tau|d\zeta|^2+\tau^{-1}|d\eta|^2\big)}$
\item
\label{blagaj}
$
Q(\omega)-\omega\cdot\nabla Q(\omega)\geqslant \delta\big(\tau |\zeta|^2+\tau^{-1}|\eta|^2\big).
$
\end{enumerate}
\end{thm}

To clarify the notation let us say that by \eqref{kupres} we mean that 
$$
\sk{-d^2Q(\omega)w}{w}\geqslant \delta\big(\tau |w_1|^2+\tau^{-1}|w_2|^2\big)
$$
for all $w=(w_1,w_2,w_3,w_4)\in\C^M\times\C^N\times\R\times\R$.

When constructing their ``scalar" Bellman function in \cite{NT}, Nazarov and Treil aimed at property \eqref{anterija} and a weaker version of \eqref{kupres}, namely ${-d^2 Q(\omega) \geqslant 2\delta |d\zeta||d\eta|}$. 
They apparently did not study anything like \eqref{blagaj}. 
It does not seem that \eqref{anterija} and \eqref{kupres} imply \eqref{blagaj}.
It was thus a 
considerable
surprise for us to see that  \eqref{blagaj} is nevertheless also true, especially since we can prove it with the same $\tau$ as in \eqref{kupres}, which is essential for our applications (see Lemma \ref{desetikvartet}).

\bigskip
\noindent
{\bf Proof.}
The inequality \eqref{anterija} is obvious. Let us first prove \eqref{kupres}.
Consider
\begin{equation*}
\label{sarmat1}
\Phi(\zeta,\eta)=|\zeta|^p+|\eta|^q+\delta\left\{
\aligned
& |\zeta|^2|\eta|^{2-q} & ; & \ \ |\zeta|^p\leqslant |\eta|^q\\
& \frac{2}{p}|\zeta|^p+\left(\frac{2}{q}-1\right)|\eta|^q
& ; &\ \ |\zeta|^p\geqslant |\eta|^q
\endaligned\right.\,.
\end{equation*}
Here (unlike before)
we think of $\zeta$ and $\eta$ as {\it real} vectors in real $l^2_{2M}, l^2_{2N}$ correspondingly. Also, $|\cdot|$ denotes the $l^2$ norm of the corresponding vector. We want to take a look at $d^2\Phi$ (the Hessian of $\Phi$, its second differential form).
To do that write $\Phi=\phi\circ U$, where $U(\zeta,\eta)=(|\zeta|,|\eta|)$ and $\phi$ is 
a function of 
two non-negative real 
variables given by
\begin{equation}
\label{sarmat2}
\phi(u,v)=u^p+v^q+\delta\left\{
\aligned
& u^2 v^{2-q} & ; & \ \ u^p\leqslant v^q\\
& \frac{2}{p}u^p+\left(\frac{2}{q}-1\right)v^q
& ; &\ \ u^p\geqslant v^q\,\ \,.
\endaligned\right.
\end{equation}

Denote by
$e_{\zeta}$ the unit vector $\zeta/|\zeta|$ and by $P_{\zeta}$ the orthogonal projection onto the orthogonal complement of $e_{\zeta}$:
$P_{\zeta}h = h - \sk{h}{e_{\zeta}}e_{\zeta}$. Then an easy direct computation gives (see \cite{NT})
$$
d|\zeta| = \sk{d\zeta}{e_{\zeta}}\hskip 20pt \text{and}\hskip 20pt d^2|\zeta| = \frac{|P_{\zeta}d\zeta|^2}{|\zeta|}\,.
$$
The same for $\eta$. 
To clarify the notation, these formul\ae\ are to be understood in the sense that if $f(a)=|a|$, then, for $a\not =0$, $df(a)h=\sk{h}{e_a}$ and $\sk{d^2f(a)h}{h}=|P_ah|^2/|a|$.

The application of the chain rule gives
\begin{equation}
\label{dPhi1}
\aligned
d^2&\Phi_{\zeta, \eta}(d\zeta, d\eta)\\
& = d^2\phi_{|\zeta|,|\eta|}(d|\zeta|, d|\eta|) + \frac{\pd \phi}{\pd u}(|\zeta|,|\eta|) d^2|\zeta| +  \frac{\pd \phi}{\pd v}(|\zeta|,|\eta|) d^2|\eta| \\
& =d^2\phi_{|\zeta|,|\eta|}(\sk{d\zeta}{e_\zeta}, \sk{d\eta}{e_{\eta}}) + \frac{\pd \phi}{\pd u}(|\zeta|,|\eta|) \frac{|P_{\zeta}d\zeta|^2}{|\zeta|} +  \frac{\pd \phi}{\pd v}(|\zeta|,|\eta|)\frac{|P_{\eta}d\eta|^2}{|\eta|}\,.
\endaligned
\end{equation}

Denote $A=|d\zeta|$. 
Since $|\sk{dh}{e_h} |^2 +|P_{h}dh|^2 =|dh|^2$ for any vectors $h,dh$,
it follows that $\sk{d\zeta}{e_\zeta}=Aa$ and $|P_{\zeta}d\zeta|^2=A^2(1-a^2)$ for some $a\in[-1,1]$. In the same way (related to $\eta,d\eta$) introduce $B,b$.
Write also $u=|\zeta|$, $v=|\eta|$.
Therefore our task is to find as good lower estimates as possible for the expression
$$
F_{u,v}(A,a,B,b):=
d^2\phi_{u,v}(Aa,Bb) + \frac{\pd \phi}{\pd u}(u,v) \frac{A^2(1-a^2)}{u} +  \frac{\pd \phi}{\pd v}(u,v)\frac{B^2(1-b^2)}{v}\,.
$$

First consider the case $u^p\geqslant v^q$. Then, by \eqref{sarmat2}, 
$$
\phi=\bigg(1+\frac2p\delta\bigg)u^p+\bigg(1+\bigg(\frac2q-1\bigg)\delta\bigg)v^q\,,
$$ 
therefore 
$$
F_{u,v}(A,a,B,b)=(p+2\delta)[1+(p-2)a^2]A^2u^{p-2}+(q+(2-q)\delta)[1-(2-q)b^2]B^2v^{q-2}\,.
$$
Note that the assumption $u^p\geqslant v^q$ implies $v^{q-2}\geqslant u^{2-p}$. Moreover, since $1+(p-2)a^2\geqslant 1$ and $1-(2-q)b^2\geqslant q-1$, 
$$
\aligned
F_{u,v}(A,a,B,b)& \geqslant(p-1)u^{p-2}A^2+(q-1)u^{2-p}B^2\\
& = \tau A^2 +\tau^{-1} B^2\,
\endaligned
$$
with $\tau=(p-1)u^{p-2}$.

Next we address the case $u^p\leqslant v^q$. This time 
$
\phi=u^p+v^q+\delta u^2v^{2-q} 
$ 
and so 
$$
F_{u,v}(A,a,B,b)=a^2U+2abV-b^2W+Z\,,
$$
where
\begin{equation*}
\aligned
U & =  p(p-2)u^{p-2}A^2\\
V  & = 2\delta(2-q)uv^{1-q}AB\\
W & = (2-q)q(\delta u^2v^{-q}+v^{q-2})B^2\\
Z & = (pu^{p-2}+2\delta v^{2-q})A^2+[(2-q)\delta u^2v^{-q}+qv^{q-2}]B^2\,.
\endaligned
\end{equation*}
These terms are all positive.

First let us fix $u,v$ and minimize $F_{u,v}(A,a,B,b)$ over all $a,b\in[-1,1]$. Since $F_{u,v}(A,-a,B,-b)=F_{u,v}(A,a,B,b)$, we may assume that $a\in[0,1]$. Then $F_{u,v}(A,a,B,-b)\leqslant F_{u,v}(A,a,B,b)$ for non-negative $b$, since $V>0$. Thus we can furthermore restrict ourselves to $b\in[-1,0]$. But this is the same as minimizing $\tilde F(a,b)=a^2U-2abV-b^2W+Z$ over $a,b\in[0,1]$.

The only stationary point of $\tilde F$ is $(0,0)$, which is obviously not the minimum. As for the boundary of $[0,1]^2$, we quickly see there are two possibilities: if $U\leqslant V$ then the minimum is attained at the point $(1,1)$ and has the value $U-2V-W+Z$, and if $U\geqslant V$ then the minimum occurs at $(V/U,1)$ and takes the value $-V^2/U-W+Z$. In both cases the minumum is minorized by the expression $Z-2V-W$.
Now $$
\aligned
Z-2V-W & = 
(pu^{p-2}+2\delta v^{2-q})A^2-4\delta(2-q)uv^{1-q}AB
\\
& \hskip 10pt +\delta[8v^{q-2}-(2-q)(q-1)u^2v^{-q}]B^2\,.
\endaligned
$$
Recall we are working under the assumption $u^p\leqslant v^q$. This gives estimates 
$u^2v^{-q}\leqslant v^{q-2}$ and $uv^{1-q}\leqslant 1$. 
Consequently
\begin{equation}
\label{polnoc}
Z-2V-W\geqslant 
\delta(2v^{2-q}A^2+7 v^{q-2}B^2)-4\delta AB
\,.
\end{equation}
For any positive $\lambda$ and $v$ we can estimate
$$
4AB\leqslant 2(\lambda v^{2-q}A^2+\lambda^{-1}v^{q-2}B^2)\,.
$$
Together with \eqref{polnoc} this implies
$$
Z-2V-W\geqslant 
\delta[2(1-\lambda)v^{2-q}A^2+(7-2\lambda^{-1}) v^{q-2}B^2]\,.
$$
By choosing, for example, $\lambda=1/2$ we get 
$$
Z-2V-W\geqslant 
\delta(v^{2-q}A^2+3v^{q-2}B^2)
\geqslant 
\delta(\tau A^2+\tau^{-1}B^2)
$$
where $\tau=v^{2-q}$.

\bigskip
Now we turn towards proving \eqref{blagaj}.
To estimate $Q(v)- v\cdot\nabla Q(v)$ from below we need to estimate $V\cdot\nabla\Phi(V) - \Phi(V)$ from below, where  $V:=(\zeta, \eta)$ is now understood as a real $(2M+2N)-$vector. We have
$$
\nabla\Phi(V)= \frac{\pd \phi}{\pd x}(|\zeta|,|\eta|)\nabla |\zeta| + \frac{\pd \phi}{\pd y}(|\zeta|,|\eta|)\nabla|\eta|\,,
$$
$$
\zeta\cdot\nabla |\zeta|=
\zeta\cdot e_\zeta=|\zeta|\hskip 15pt \text{and}\hskip 15pt \eta\cdot\nabla |\eta|=
\eta\cdot e_\eta=|\eta|\,.
$$
Combining this and writing again $u=|\zeta|$, $v=|\eta|$, we get
$$
V\cdot\nabla\Phi(V)- \Phi(V)= \frac{\pd \phi}{\pd u}(u,v) u+\frac{\pd \phi}{\pd v}(u,v)v -\phi(u,v)\,.
$$
Denote this expression by $\Lambda(u,v)$. 

By now, having finished the proof of part \eqref{kupres}, we already have candidates for $\tau$ we must work with. 

Suppose $u^p\geqslant v^q$. Then a direct calculation shows that 
$$
\aligned
\Lambda(u,v) & = \Big(p-1+\frac{q-1}4\Big) u^p +(q-1)\Big(1+\frac{(2-q)(q-1)}8\Big)v^q\\
& \geqslant (p-1)u^{p-2}u^2 +(q-1)v^{q-2}v^2\,.
\endaligned
$$
But $u^p\geqslant v^q$ implies  $v^{q-2}\geqslant u^{2-p}$, which proves $\Lambda(u,v)\geqslant \tau u^2 +\tau^{-1}v^2$ with the same $\tau$ as in the corresponding case of \eqref{kupres}, namely, $\tau=(p-1)u^{p-2}$.

Finally, suppose $u^p\leqslant v^q$. Then 
$$
\Lambda(u,v)  = (p-1) u^p +(q-1)v^q+(3-q)\delta u^2v^{2-q}\,.
$$
The first term on the right we simply drop out. Since $q-1>\delta$ and $3-q\geqslant 1$ we get 
$$
\Lambda(u,v)\geqslant \delta(v^{q-2}v^2+v^{2-q}u^2)\,
$$
so the required inequality is proven with $\tau=v^{2-q}$, exactly as in the proof of \eqref{kupres}. 
\qed

\bigskip
\noindent
{\bf Remark.} 
In the proof we saw that when $|\zeta|^p\geqslant |\eta|^q$ 
we get
the properties \eqref{kupres} and \eqref{blagaj} without $\delta$.

\bigskip

We also need to know how the gradient of $Q$ behaves. Calculations carried out on the basis of \eqref{sarmat} give estimates
\begin{equation}
\label{brcko}
 \Big|\frac{\pd Q}{\pd \zeta}\Big|
 \leqslant
C(p)\max\{|\zeta|^{p-1},|\eta|\}
\hskip 15pt\text{and}\hskip 15pt
\Big|{\pd Q\over \pd \eta}\Big|
\leqslant
C|\eta|^{q-1}\,.
\end{equation}
Here of course
$$
\frac{\pd Q}{\pd \zeta}=\bigg(\frac{\pd Q}{\pd \zeta_1},\hdots,\frac{\pd Q}{\pd \zeta_M}\bigg)
\hskip 15pt\text{and}\hskip 15pt
\frac{\pd Q}{\pd \eta}=\bigg(\frac{\pd Q}{\pd \eta_1},\hdots,\frac{\pd Q}{\pd \eta_N}\bigg)\,.
$$
Same estimates apply to the $\bar\pd-$derivatives of $\zeta$ and $\eta$, for $Q$ is a real-valued function.

\section{Bilinear embedding}

This section is devoted to proving Theorem \ref{bilinher}.
At this point it might be worth explaining the origins of the term ``bilinear embedding" which we use throughout the article.

Let $\phi$ be a $C^1$ complex-valued function, defined on $\R^{n}\times (0,\infty)$.
Write
$$
\nabla_*\phi(x,t)=(\nabla \phi(x,t),x\,\phi(x,t))\in\C^{2n+1}\,.
$$
Then the statement of Theorem \ref{bilinher} 
implies
the following one:

\medskip
\noindent
{\it The pairing $\Psi$, given by
$$
\Psi(f,g)(x,t)=\nabla_* P_t f(x) \cdot \nabla_* P_t g(x)\,,
$$
defines a bounded bilinear
mapping $L^p\times L^q \rightarrow L^1(\R^n\times
(0,\infty),\, dx\,t\,dt)$ for all $p\in (1, \infty)$.
Its norm 
is controlled by $C(p^*-1)$. 
}


\subsection{Structure of the proofs}

\label{plan}

We already emphasized that one of the main features of this presentation is the uniformity of the proofs regardless of the semigroup we work with. Having already introduced the Bellman function and its properties, we now wish to illustrate our strategy a bit further. 

 \medskip
Given test functions $f,g$ on $\R^n$, we want to define
\begin{equation}
\label{detinec}
v(x,t):= (P_t f(x), P_t g(x), P_t |f|^p(x), P_t |g|^q(x))\,
\end{equation}
and furthermore $b:=Q\circ v$, that is,
\begin{equation}
\label{mersada}
b(x,t) := Q(P_t f(x), P_t g(x) , P_t |f|^p(x), P_t |g|^q(x))\,.
\end{equation}
For that purpose we have to check two things. One is that 
$P_t|f|^p$ is well defined.
The other is that $v(x,t)\in\Omega$. This is true for most of the natural semigroup extensions in which one can express $P_t\varphi(x)$ as an integral of $\varphi$ against some 
finite
measure depending on $(x,t)$
(see \cite{DV} for classical and Gaussian case and \cite{T} for Hermite operator).
This permits the inequality $|P_t\varphi|^p\leqslant P_t|\varphi|^p$, which is exactly what we need. Having explicit formulas at our disposal will settle the questions of well-posedness of $b$ for our purposes. In general, though, it is not known to us for what class of $\sqrt{L}$-generated extensions this holds.

The main two steps in the proofs of the embedding theorems are always the same:
\begin{itemize}
\item
Consider the operator
$$
\label{borjomi}
L'=\frac{\pd^2}{\pd t^2}-L\,.
$$
It can be regarded as an extension of $-L$ in the upper half-space.
This extension is the ``right one", in the sense that $L'\widetilde\varphi=0$. Namely, this enables us to express $L'b$ in terms of the Hessian of $Q$, and thereupon everything is set for applying the concavity properties \eqref{kupres} and \eqref{blagaj} of $Q$.

\item
Our aim is to estimate the integral
\begin{equation}
\label{hunedoara}
-\int L' b(x,t)\,dx\,t\, dt\,
\end{equation}
from below and above.
The size property \eqref{anterija} of $Q$ makes up for the upper estimate of the integral above, whereas \eqref{kupres} and \eqref{blagaj} provide the estimates from below. The expressions which appear in the lower and the upper estimate of \eqref{hunedoara} are exactly those from the embedding theorem. That is, a more complete formulation of Theorem \ref{bilinher} 
would incorporate \eqref{hunedoara} as the middle term in the inequalities.

\end{itemize}

We carry out the plan described above. 
The particularity of the Hermite case is hidden in the fact that the said operator has a potential, $|x|^2$, which has to be reckoned with. This means that, in contrast with 
the situations studied in \cite{DV}, the formula for $L'b$ contains not only the scalar products involving the Hessian of $Q$, but also some other terms. However, the property \eqref{blagaj} of $Q$ gives control exactly over these, newly arisen terms.

The proof of the Lemma below makes this statement more transparent.

\subsection{Estimate of the integral \eqref{hunedoara} from below}
 
The lower estimates of the integral \eqref{hunedoara} will trivially follow from the lower {\sl pointwise} estimates of $L'b(x,t)$ which we present next. Recall the notation $\delta=q(q-1)/8$.

\begin{lema}
\label{desetikvartet}
There is an absolute $C>0$ such that for all test functions $f,g$ and all $p\geqslant 2$ we have 
\begin{equation}
\label{hasiba}
- L'b(x,t) \geqslant 
\delta\nor{P_tf(x)}_{*}\nor{P_tg(x)}_{*}\,.
\end{equation}
\end{lema}
\dok
By applying the chain rule  we get
\begin{equation}
\label{herculian}
- L'b(x,t) = \sum_{j=0}^n \Sk{-d^2 Q(v_0) \frac{\pd v}{\pd x_j}(x,t)}{\frac{\pd v}{\pd x_j}(x,t)} 
+ |x|^2 [Q(v_0) - v_0\cdot \nabla Q(v_0)]\,.
\end{equation}
Here we wrote $v_0=v(x,t)$ and when $j=0$ we meant the differentiation in $t$.
Combining \eqref{herculian} and Theorem \ref{bejaz}  we find $\tau=\tau(x,t)>0$ such that
$$
\aligned
- L'b(x,t) 
\geqslant & \delta\tau\Big(\sum_{j=0}^n \Big|\frac{\pd}{\pd x_j} P_t f(x)\Big|^2 +  |x|^2 |P_t f(x)|^2\Big)
\\ 
& 
+ \delta\tau^{-1}\Big(\sum_{j=0}^n \Big|\frac{\pd}{\pd x_j} P_tg(x)\Big|^2 + |x|^2 |P_tg(x)|^2\Big)
\,.
\endaligned
$$
Now the inequality between the arithmetic and the geometric mean gives \eqref{hasiba}. \qed

\bigskip
\noindent
{\bf Remark.} 
\label{busno1}
Notice that if  $|x|^2$ is  replaced by some arbitrary non-negative measurable function $V(x)$ in the definition of $L$ (and consequently in $L'$, $P_t$ and $\nor{\cdot}_{*}$), the statement from Lemma \ref{desetikvartet} does not change at all.

\subsection{Estimates of the semigroup kernels}

The aim of this section is to compile some of the known estimates for the integral kernels of the heat and Poisson semigroup, respectively, 
and to adapt them to suit our purposes. The estimates we have in mind are pointwise estimates of the kernel and its derivatives. 
We will need them in the continuation of the proof of Theorem \ref{bilinher}, that is, when giving upper estimates of the integral $-\int L'b$ introduced in \eqref{hunedoara}.

\subsubsection{Estimates of the Hermite heat kernel}

From \cite[4.1.2]{T} we have that
\begin{equation}
\label{8-2}
e^{-tL}\f(x)={1\over (2\pi)^{n/2}}\int_{\R^n}K_t(x,y)\f(y)\,dy\,,
\end{equation}
where \cite[4.1.3]{T}
\begin{equation}
\label{harlech}
K_t(x,y)={1\over (\sinh 2t)^{n/2}}\,\exp\Big(-\frac{|x|^2+|y|^2}{2}\coth2t+{\sk{x}{y}\over \sinh 2t}\Big)
\,. 
\end{equation}
Direct calculation shows that
\begin{equation}
\label{heat}
K_t(x,y) \leqslant (2t)^{-\frac{n}{2}}\,{e^{-\frac{|x-y|^2}{4t}}} \,,
\end{equation}
so the Hermite heat kernel is majorized by the ordinary heat kernel.

By Lemma 4.3.1.(i) and Lemma 4.3.2.(ii) from \cite{T} there exist positive constants $C,a $ 
not depending on $x,y,t$ such that for $j\in\{1,\hdots,n\}$ we have
\begin{equation}
\label{ortlinde}
\bigg|{\pd \over \pd x_j}K_{t}(x,y)\bigg|
\leqslant C
t^{-{n+1\over 2}}\,e^{-{a\over t}|x-y|^2} 
\,,
\end{equation}
whereas 
Lemma 4.1.1 (i) from the same source gives 
\begin{equation}
 \label{heatderiv}
\bigg|\frac{\pd}{\pd t} K_t(x,y)\bigg| \leqslant Ct^{-\frac{n}{2}-1}\, e^{-\frac{a}{t}|x-y|^2}\,.
 \end{equation}
These estimates are valid for all $t>0$.

\subsubsection{Estimates of the Hermite Poisson kernel}

\label{goggles}

Most of the inequalities just encountered can be transferred to the Poisson semigroup by the following well-known subordination principle. 

Denote
\begin{equation*}
\label{gurck}
d\mu(s)={1\over \sqrt\pi}\,e^{-s}s^{-1/2}\,ds\,.
\end{equation*}
This is a probability measure on $(0,\infty)$. The integral equation
$$
e^{-\alpha}=\int_0^\infty e^{-{\alpha^2\over 4s}}d\mu(s)
$$
gives rise to the subordination formula
\begin{equation}
\label{assai}
P_t\f(x)=\int_0^\infty e^{-{t^2\over 4s}L}\f(x)\,d\mu(s) \,.
\end{equation}

By the same symbol as the extension, namely $P_t$,
we also denote the Hermite Poisson {\sl kernel}. 
From \eqref{assai} and \eqref{8-2} 
we get
\begin{equation}
\label{armstrong}
P_t(x,y) = \frac{1}{(2\pi)^{n/2}}\int_0^{\infty} K_{\frac{t^2}{4s}}(x,y)\, d\mu(s)\,. 
\end{equation}
Hence, by \eqref{heat},
\begin{equation}
\label{lance}
\aligned
P_t(x,y) 
\leqslant \frac{\Gamma\big(\frac{n+1}2\big)}{\pi^{\frac{n+1}2}}\cdot\frac{t}{(|x-y|^2 + t^2)^{\frac{n+1}{2}}}\,.
\endaligned
\end{equation}
In other words our rough estimate shows that $P_t(x,y)$ is majorized by the classical Poisson kernel in $\R^n$.

In order to give estimates of the derivatives of $P_t$, combine 
\eqref{armstrong} and \eqref{ortlinde} to deduce
\begin{equation}
\label{hvar}
\Big|{\pd \over \pd x_j}P_t(x,y)\Big|
\leqslant 
\frac{C t}{(|x-y|^2 + t^2)^{\frac{n+2}{2}}}\,.
\end{equation}
As for the estimate of the derivative in $t$, 
fix $(x,y)$ and write $\psi(u):=K_u(x,y)$. Then
\begin{equation}
\label{zapiski}
{\pd \over \pd t}[K_{\lambda t^2}(x,y)]=2\lambda t \psi'(\lambda t^2)\,.
\end{equation}
We can calculate $\psi'$ from \eqref{harlech}:
\begin{equation}
\label{sejo}
\psi'(u)=\psi(u)\left[{|x|^2+|y|^2-2\sk{x}{y}\cosh 2u \over \sinh^2 2u} -n\coth 2u\right]\,.
\end{equation}
The estimate for $\psi'$, however, comes from \eqref{heatderiv}. We apply it together with \eqref{armstrong} and \eqref{zapiski}. The result is 
\begin{equation}
\label{brickroad}
\Big|{\pd \over \pd t}P_t(x,y)\Big|
\leqslant 
\frac{C}{(|x-y|^2 + t^2)^{\frac{n+1}{2}}}\,.
\end{equation}

\bigskip
\noindent
{\bf Remark.}
It follows from \eqref{heat} that $(2\pi)^{-n/2}K_t(x,y)\,dy$ is, for each $x\in\R^n$, a positive sub-probability measure on $\R^n$ (later on, see \eqref{abc}, we precisely calculate its mass). 
Jensen's inequality implies from here the pointwise estimate $|e^{-tL}\f(x)|^p\leqslant e^{-tL}|\f|^p(x)$ whenever $p\geqslant 1$. By \eqref{assai}, the same is true for the Poisson semigroup, i.e. $|P_t\f(x)|^p\leqslant P_t|\f|^p(x)$. 
This in retrospect justifies our definition of function $b$, i.e. vector $v(x,t)$ really maps into the domain of the Bellman function $Q$; see the
discussion ensuing the definitions
\eqref{detinec}, \eqref{mersada}.

\subsection{Integration by parts}
\label{ibp}

In this section our goal is to extract the ``noncontributing" part of the integral \eqref{hunedoara}.
\begin{lema}
\label{fazila}
Let $f$ and $g$ belong to 
 $C_c^\infty$ and let $b$ be as in \eqref{mersada}.
For all $t>0$,
\begin{equation}
\label{danger2}
-\int_{\R^n}  L' b(x,t)\, dx= \int_{\R^n}\Big(-\frac{\pd^2b}{\pd t^2}(x,t) + |x|^2 b(x,t)\Big)\,dx\,.
\end{equation}
\end{lema}
\dok
Write
\begin{equation*}
\label{scherzo}
I:= 
\int_{\R^n} \Delta b(x,t)
\, dx
\,.
\end{equation*}
Then \eqref{danger2} translates
into showing that $I=0$. By symmetry it suffices to do that for $\pd^2/\pd x_1^2$ in place of $\Delta$.

Take $M>0$ and denote $R_M=[-M,M]^n$. If  $x=(x_1,\hdots,x_n)$, denote temporarily $x'=(x_2,\hdots,x_n)$, so that we can write $x=(x_1,x')$. Then  
$$
\aligned
\int_{R_M}\frac{\partial^2b}{\partial x_1^2}& (x,t)\,dx_1\hdots dx_n\\
& =\underbrace{\int_{-M}^M\hdots\int_{-M}^M}_{n-1}
\Big[\frac{\partial b}{\partial x_1}(M,x',t)-\frac{\partial b}{\partial x_1}(-M,x',t)\Big]
\,dx_2\hdots dx_n\,.
\endaligned
$$
Thus
\begin{equation}
\label{gerhilde}
\bigg|\int_{R_M}\frac{\partial^2b}{\partial x_1^2} (x,t)\,dx_1\hdots dx_n\bigg|
\leqslant2S\,(2M)^{n-1}\,,
\end{equation}
where 
$$
S=\sup_{x\in\partial R_M}\bigg|\frac{\partial b}{\partial x_1} (x,t)\bigg|\,.
$$

We want to estimate $S$. 
As in \cite{DV},
\begin{equation*}
{\pd b\over \pd x_1}(x,t)=\Sk{\nabla Q(v_0)}{{\pd v\over \pd x_1}(x,t)}\,,
\end{equation*}
where $v_0=v(x,t)$. This means
%
%
\begin{equation}
\label{marcia}
\aligned
{\pd b\over \pd x_1}(x,t)
=&\Sk{{\pd Q\over \pd \zeta}(v_0)}{{\pd \over \pd x_1}P_tf(x)}_{\C^M}
+\Sk{{\pd Q\over \pd \eta}(v_0)}{{\pd \over \pd x_1}P_tg(x)}_{\C^N}\\
+&\Sk{{\pd Q\over \pd \bar\zeta}(v_0)}{{\pd \over \pd x_1}\overline{P_tf(x)}}_{\C^M}
+\Sk{{\pd Q\over \pd \bar\eta}(v_0)}{{\pd \over \pd x_1}\overline{P_tg(x)}}_{\C^N}\\
&+{\pd Q\over \pd Z}(v_0){\pd \over \pd x_1}P_t|f|^p(x)
+{\pd Q\over \pd H}(v_0){\pd \over \pd x_1}P_t|g|^q(x)\,.
\endaligned
\end{equation}
%
Since in our case $\zeta=P_tf(x)$ and $\eta=P_tg(x)$, \eqref{brcko} implies
\begin{equation}
\label{siegrune}
\aligned
\bigg|\frac{\partial b}{\partial x_1} (x,t)\bigg|
\leqslant 
& C\bigg(\big(|P_tf(x)|^{p-1}+|P_tg(x)|\big)\Big|{\pd\over\pd x_1}P_tf(x)\Big| \\ 
& +|P_tg(x)|^{q-1}\Big|{\pd\over\pd x_1}P_tg(x)\Big|
+\Big|{\pd\over\pd x_1}P_t|f|^p(x)\Big|+\Big|{\pd\over\pd x_1}P_t|g|^q(x)\Big|\bigg)\,.
\endaligned
\end{equation}
Hence the estimation of $S$ is reduced to estimating $P_t\varphi(x)$ and ${\pd\over\pd x_1}P_t\varphi(x)$
with $f,g,|f|^p,|g|^q$ in place of $\varphi$. In order to do that we 
recall the estimates from section \ref{goggles}.

Let the radius $A>0$ be such that the ball $B(0,A)$ contains $\text{supp }\f$.
If
$y\in\text{supp }\f\subseteq B(0,A)$ and $|x|\geqslant 2A$ we get $|x-y|\geqslant |x|/2$.
Then 
$
P_t\f(x)=\int_{B(0,A)}P_t(x,y)\f(y)\,dy
$ 
and so for $|x|\geqslant 2A$ the inequality \eqref{lance} implies  
\begin{equation}
\label{grimgerde}
|P_t\f(x)|\leqslant \frac{C(n)t\nor{\f}_1}{|x|^{n+1}}\,.
\end{equation}

Now we turn to the estimate of ${\pd\over\pd x_1}P_t\f(x)$. 
First,
\begin{equation*}
\label{helmwige}
\Big|{\pd \over \pd x_1}P_t\f(x)\Big|
\leqslant\int_{\R^n}\Big|{\pd \over \pd x_1}P_t(x,y)\Big|\,|\f(y)|\,dy
\,.
\end{equation*}
The inequality \eqref{hvar} implies
\begin{equation}
\label{rossweisse}
\Big|{\pd \over \pd x_1}P_t\f(x)\Big|
\leqslant
\frac{C(n)t\nor{\f}_1}{\big(|x|^2+t^2\big)^{\frac{n+2}2}}\,
\end{equation}
for sufficiently large $x$ as specified above.

Note that $x\in\pd R_M$ implies $|x|\geqslant M$. 
Now a combination of \eqref{gerhilde}, \eqref{siegrune}, \eqref{grimgerde} and \eqref{rossweisse} shows that 
$$
\lim_{M\rightarrow\infty}\int_{R_M}\frac{\partial^2b}{\partial x_1^2} (x,t)\,dx_1\hdots dx_n
=0\,.
$$
This proves \eqref{danger2} and thus Lemma \ref{fazila}. \qed

\subsection{Estimate of the integral \eqref{hunedoara} from above}
\label{pritvoriseplaninata}

Here we treat a consequence of Lemma \ref{fazila} which consists of showing that the expression in \eqref{danger2}, and for that matter the integral \eqref{hunedoara}, are bounded by $C(\nor{f}_p^p+\nor{g}_q^q)$.

\begin{prop}
\label{manco}
For all $p>1$,
$$
-\int_0^\infty\int_{\R^n}  L' b(x,t)\,t\, dx\, dt
\leqslant
6(\nor{f}_p^p+\nor{g}_q^q)\,.
$$
\end{prop}
\dok
It clearly suffices to consider the case $p\geqslant 2$.
By Lemma \ref{fazila} we are done once we prove 
\begin{equation}
\label{dety}
\int |x|^2b(x,t)\,dx\, t\, dt \leqslant 4(\nor{f}_p^p+\nor{g}_q^q)
\end{equation}
and
\begin{equation}
\label{dittie}
-\int\frac{\pd^2 b}{\pd t^2}(x,t)\,dx \,t\,dt\leqslant 2(\nor{f}_p^p+\nor{g}_q^q)\,.
\end{equation}

\subsubsection{
\label{x^2}
Proof of \eqref{dety}.
}

It follows from \eqref{anterija} on page \pageref{anterija} that
\begin{equation}
\label{certina}
\int |x|^2b(x,t)\,dx\, t\, dt\leqslant 2\int |x|^2 (P_t|f|^p(x)+P_t|g|^q(x))\,dx\, t\, dt\,.
\end{equation}
The combination of \eqref{assai}, \eqref{8-2} 
gives
$$
\aligned
\int_0^\infty\int_{\R^n} & |x|^2 P_t|f|^p(x)\,dx\, t\, dt\\
& =\int_0^\infty\int_{\R^n} |x|^2 \int_0^\infty e^{-{t^2\over 4s}L}|f|^p(x)\,d\mu(s)\,dx\, t\, dt\\
& =\int_0^\infty\int_{\R^n} |x|^2 \int_0^\infty {1\over (2\pi)^{n/2}}\int_{\R^n}
K_{t^2\over 4s}(x,y)|f(y)|^p\,dy\,d\mu(s)\,dx\, t\, dt
\,.
\endaligned
$$
It should suffice for our purpose to know that the integrals
\begin{equation}
\label{dunvegan}
{1\over (2\pi)^{n/2}}\int_0^\infty\int_{\R^n} \int_0^\infty
|x|^2K_{t^2\over 4s}(x,y)\,
\,d\mu(s)\,dx\, t\, dt
\end{equation}
are bounded uniformly in $y$ and $n$. 
By \eqref{harlech}, the expression in \eqref{dunvegan} equals
$$
\aligned
{1\over (2\pi)^{n/2}}\int_0^\infty&\int_{\R^n} \int_0^\infty
{|x|^2\over (\sinh {t^2\over 2s})^{n/2}}\,\cdot\\
& \ \hskip 10pt\cdot \exp\Big(-\frac{|x|^2+|y|^2}{2}\coth{t^2\over 2s}+{\sk{x}{y}\over \sinh{t^2\over 2s}}\Big)
\,d\mu(s)\,dx\, t\, dt\,.
\endaligned
$$

First we integrate in $x$. Write temporarily $\alpha=\coth{t^2\over 2s}$, $\beta=\sinh{t^2\over 2s}$. Consider
$$
\hskip -60pt
\int_{\R^n}
|x|^2\, \exp\Big(-\frac{|x|^2+|y|^2}{2}\alpha+{\sk{x}{y}\over \beta}\Big)
\,dx\,
$$
\begin{equation}
\label{beaumaris}
\hskip 70pt=\int_{\R^n}
\sum_{j=1}^nx_j^2\, \exp\Big(-\sum_{k=1}^n\Big(\frac{x_k^2+y_k^2}{2}\alpha-{x_ky_k\over \beta}\Big)\Big)
\,dx\,.
\end{equation}
Note that
$$
\frac\alpha2(x_k^2+y_k^2)-{x_ky_k\over \beta}
=\frac12\Big(\sqrt\alpha x_k-{y_k\over\sqrt\alpha\beta}\Big)^2+{y_k^2\over 2\alpha}\,,
$$
hence we can continue \eqref{beaumaris} as
\begin{equation}
\label{kisimul}
\aligned
e^{-{|y|^2\over 2\alpha}}
\sum_{j=1}^n\int_{\R^n}x_j^2\,\prod_{k=1}^n e^{-\frac12\big(\sqrt\alpha\,x_k-{y_k\over\sqrt\alpha\,\beta}\big)^2}
\,dx_1\hdots dx_n\,.
\endaligned
\end{equation}
Now use that
$$
\int_{\R}
x^2\,e^{-\frac12\big(\sqrt\alpha\,x-{y\over\sqrt\alpha\,\beta}\big)^2}
\,dx = \sqrt{2\pi\over\alpha} \Big( {y^2\over \alpha^2\beta^2}+\frac1\alpha
\Big)
$$
and
$$
\int_{\R}
\,e^{-\frac12\big(\sqrt\alpha\, x-{y\over\sqrt\alpha\,\beta}\big)^2}
\,dx=\sqrt{2\pi\over\alpha}\,,
$$
which simplifies \eqref{kisimul} to
$$
e^{-{|y|^2\over 2\alpha}}\Big({2\pi\over\alpha}\Big)^{n/2}\Big( {|y|^2\over \alpha^2\beta^2}+\frac n\alpha
\Big)\,.
$$
Therefore we proved that the integral in \eqref{dunvegan} is equal to
$$
{1\over (2\pi)^{n/2}}\int_0^\infty \int_0^\infty
{1\over \beta^{n/2}}\,
e^{-{|y|^2\over 2\alpha}}\Big({2\pi\over\alpha}\Big)^{n/2}\bigg( {|y|^2\over \alpha^2\beta^2}+\frac n\alpha
\bigg)
\,d\mu(s)\, t\, dt
$$
or
\begin{equation}
\label{doune}
\int_0^\infty \int_0^\infty
{e^{-{|y|^2\over 2}\tanh{t^2\over 2s}}\over \cosh^{1+n/2}{t^2\over 2s}}\,
\bigg( {|y|^2\over \cosh{t^2\over 2s}}+n\sinh{t^2\over 2s}
\bigg)
\,d\mu(s)\, t\, dt\,.
\end{equation}
Introduce a new variable $u=t^2/2$ and write \eqref{doune} as 
$I_1+I_2$, where
$$
\aligned
I_1= &
\int_0^\infty \int_0^\infty
|y|^2e^{-{|y|^2\over 2}\tanh(u/s)}\cosh^{-2-n/2}(u/s)\,
\,d\mu(s)\,
du\\
I_2= &
\int_0^\infty \int_0^\infty
e^{-{|y|^2\over 2}\tanh(u/s)}n\cosh^{-1-n/2}(u/s)\,
\sinh(u/s)
\,d\mu(s)\,
du\,.
\endaligned
$$
Let us estimate $I_1$ first. Obviously
$$
I_1\leqslant
\int_0^\infty \int_0^\infty
|y|^2e^{-{|y|^2\over 2}\tanh(u/s)}\cosh^{-2}(u/s)
\, du\,d\mu(s)\,.
$$
In the inner integral introduce a new variable $w=e^{-{|y|^2\over 2}\tanh(u/s)}$, from where we continue with
$$
\aligned
=\int_0^\infty 2s \int_{e^{-{|y|^2\over 2}}}^1
dw
\,d\mu(s)
&=2\Big(1-e^{-{|y|^2\over 2}}\Big){1\over\sqrt\pi}\int_0^\infty e^{-s}s^{1/2}
\,ds\\
& =1-e^{-{|y|^2\over 2}}\,.
\endaligned
$$

As for $I_2$, we have
$$
I_2\leqslant
n\int_0^\infty \int_0^\infty
\cosh^{-1-n/2}(u/s)\,
\sinh(u/s)\,
du\,d\mu(s)\,.
$$
This time we take $w=\cosh^{-n/2}(u/s)$, which gives
$$
\int_0^\infty
2s\int_0^1
dw
\,d\mu(s)=1\,.
$$
So we showed that
$$
I_1+I_2\leqslant 2-e^{-{|y|^2\over 2}}<2
\,.
$$
Consequently
$$
\int |x|^2 (P_t|f|^p(x)+P_t|g|^q(x))\,dx\, t\, dt
\leqslant 2(\nor{f}_p^p+\nor{g}_q^q)\,,
$$
which, in view of \eqref{certina}, implies \eqref{dety}.

\bigskip
\noindent
{\bf Remark.}
This time we 
had to work with the exact formula for the kernels $K_t$; the approximation \eqref{heat} seems already a step to far, one where too much information is lost already.

\subsubsection{Proof of \eqref{dittie}}

It remains to prove 
a similar
bound for the 
integral 
in \eqref{dittie}. 
We follow the considerations from \cite{DV}.

Integration by parts gives
\begin{equation*}
\label{velvia}
\aligned
& -\limsup_{\omega_1\rightarrow 0\atop \omega_2\rightarrow \infty}  \int_{\R^n}\int_{\omega_1}^{\omega_2}\frac{\pd^2 b}{\pd t^2}(x,t) \,t\,dt\,dx
=\, \underbrace{-\limsup_{t\rightarrow\infty}\int_{\R^n}t\,\frac{\pd b}{\pd t}(x,t)\,dx}_{\rm I}\\
&
+\underbrace{\liminf_{t\rightarrow 0}\int_{\R^n}t\,\frac{\pd b}{\pd t}(x,t)\,dx}_{\rm II}
+\underbrace{\liminf_{t\rightarrow \infty}\int_{\R^n}b(x,t)\,dx}_{\rm III}
-\underbrace{\limsup_{t\rightarrow 0}\int_{\R^n}b(x,t)\,dx}_{\rm IV}
\,.
\endaligned
\end{equation*}


\bigskip
Let us start with $\rim{III}$. 
%
%
By \eqref{mersada} and the property \eqref{anterija} of $Q$ we have
$$
\aligned
\liminf_{t\rightarrow \infty}\int_{\R^n}b(x,t)\,dx & \leqslant 2\liminf_{t\rightarrow \infty} \int_{\R^n}(P_t|f|^p(x)+P_t|g|^q(x))\,dx\\
& \leqslant 2\int_{\R^n}(|f|^p(x)+|g|^q(x))\,dx\,,
\endaligned
$$
the second inequality following from the contractivity property of $P_t$.

Since function $Q$ (and therefore $b$) is positive, then so is the term $\rim{IV}$, therefore we can skip it from all the estimates from above.

\bigskip
We are left with $\rim{II}$ and $\rim{I}$. 
We want to show that they 
can be neglected. More precisely, we will prove that $\rim{II}\leqslant 0$ and 
$\rim{I}=0$. 

In the estimates of $\rim{II}$ and $\rim{I}$ we will continuously be applying and referring to 
the decomposition of 
${\pd b\over \pd t}(x,t)$ as in \eqref{marcia}, just with $t$ in place of $x_1$. This essentially gives four terms (nominally there are six of them, but $\pd_\zeta$ and $\pd_{\bar\zeta}$ derivatives are handled in the same fashion; likewise for $\pd_\eta$ and $\pd_{\bar\eta}$). 

First we treat the terms in $\rim{II}$ and $\rim{I}$ corresponding,
in the sense of \eqref{marcia} explained above,
to 
the partial derivatives with respect to $Z$ and $H$. These derivatives 
are identically equal to 2, hence for the corresponding terms the 
estimate
reduces 
to showing that, for $\varphi=|f|^p+|g|^q$,
\begin{equation}
\label{ECCisDEsD}
\int_{\R^n}{\pd \over \pd t}P_t\varphi(x)\,dx\leqslant 0\,,
\end{equation}
which will imply that $\rim{II}\leqslant 0$ (more precisely, its part associated with the $Z$- and $H$- derivatives, see \eqref{marcia} on p. \pageref{marcia}), 
and 
\begin{equation}
\label{pozno}
\lim_{t\rightarrow\infty}t\int_{\R^n}{\pd \over \pd t}P_t\varphi(x)\,dx=0\,,
\end{equation}
which will in turn imply that the ``$Z$- and $H$- part" of $\rim{I}$ vanishes.
Proving that is not as straightforward as in \cite{DV}, since we cannot use the scalar product (i.e. duality) argument. Instead, we resort once again to the explicit formulas for the kernels.
We start with the calculation of the integral appearing in \eqref{ECCisDEsD} and \eqref{pozno}.


By \eqref{assai} and \eqref{8-2},
\begin{equation}
\label{mujohrnjica}
\int_{\R^n}{\pd \over \pd t}P_t\f(x)\,dx={1\over (2\pi)^{n/2}}\int_{\R^n}\int_0^\infty\int_{\R^n}{\pd \over \pd t}K_{t^2\over 4s}(x,y)\f(y)\,dy
\,d\mu(s)\,dx\,.
\end{equation}
First we integrate with respect to $x$, i.e. we compute
\begin{equation}
\label{schmied}
\int_{\R^n}{\pd \over \pd t}K_{t^2\over 4s}(x,y)\,dx\,.
\end{equation}
The formul\ae\ \eqref{zapiski} and \eqref{sejo} show that the integral in \eqref{schmied} can be written as
$$
\frac{t}{2s\alpha^2}
\int_{\R^n}K_{t^2\over 4s}(x,y)(|x|^2+|y|^2-2\sk{x}{y}\beta)\,dx
-
\frac{nt\beta}{2s\alpha}
\int_{\R^n}K_{t^2\over 4s}(x,y)\,dx\,,
$$
where now $\displaystyle{\alpha=\sinh{t^2\over 2s}}$ and $\displaystyle{\beta=\cosh{t^2\over 2s}}$.

A computation shows this is the same as
$$
\frac{t}{2s\alpha^2}e^{-{\alpha|y|^2\over 2\beta}}\left({2\pi\over\beta}\right)^{n/2}\Big[-{\alpha^2\over\beta^2}|y|^2+{n\alpha\over\beta}\Big]
-
\frac{nt\beta}{2s\alpha}e^{-{\alpha|y|^2\over 2\beta}}\left({2\pi\over\beta}\right)^{n/2}\,,
$$
thus
\begin{equation}
\label{torto}
\int_{\R^n}\frac\pd{\pd t}K_{t^2\over 4s}(x,y)\,dx
=
-\frac{t}{2s}\,e^{-\frac{\alpha|y|^2}{2\beta}}\left(\frac{2\pi}{\beta}\right)^{n/2}
\left[
{|y|^2\over\beta^2}
+
\frac{n\alpha}{\beta}
\right]\,.
\end{equation}
So we have, by \eqref{mujohrnjica} and \eqref{torto},
\begin{equation}
\label{kura}
\int_{\R^n}{\pd \over \pd t}P_t\f(x)\,dx
=-\int_0^\infty{t\over 2s\beta^{n/2}}\int_{\R^n}
e^{-{\alpha|y|^2\over 2\beta}}\left[
{|y|^2\over\beta^2}
+
\frac{n\alpha}{\beta}
\right]
\f(y)\,dy
\,d\mu(s)\,.
\end{equation}
Note that the integrand is almost identical to the one in \eqref{doune}. 
More significantly, the above expression is non-positive, because $\varphi\geqslant 0$. 
This immediately implies \eqref{ECCisDEsD}.

As for 
\eqref{pozno},
first write $\f(y)\leqslant\nor{\f}_\infty$ (of course we may assume that $\f$ is bounded). One can pass to polar coordinates in $\R^n$ and explicitly calculate 
$$
\int_{\R^n}
e^{-{\alpha|y|^2\over 2\beta}}
\left[{|y|^2\over\beta^2}+\frac{n\alpha}{\beta}\right]
\,dy
=
(2\pi)^{n/2}n\bigg(\frac{\beta}{\alpha}\bigg)^{1+n/2}\,.
$$
Hence \eqref{kura} gives
\begin{equation*}
0\leqslant-t\int_{\R^n}{\pd \over \pd t}P_t\f(x)\,dx
\leqslant\nor{\f}_\infty
(2\pi)^{n/2}n
\int_0^\infty
{t^2\over 2s}\cdot \frac{\beta}{\alpha^{\frac n2+1}}
\,d\mu(s)\,.
\end{equation*}
Now 
\begin{equation*}
\label{teskoovozivot}
\frac{\beta}{\alpha^{\frac n2+1}}
=\frac{1+\alpha^2}{\beta\alpha^{\frac n2+1}}
=\frac{1}{\beta\alpha^{\frac n2+1}}+\frac{1}{\beta\alpha^{\frac n2-1}}
\,,
\end{equation*}
therefore
\begin{equation}
\label{drugadalmatinska}
\int_0^\infty
{t^2\over 2s}\cdot \frac{\beta}{\alpha^{\frac n2+1}}
\,d\mu(s)
\leqslant
\int_0^\infty\bigg(\frac{2s}{t^2}\bigg)^{\frac n2+1}\,d\mu(s)
+
C\int_0^\infty \bigg(\frac{2s}{t^2}\bigg)^{\frac n2}\,d\mu(s)\,,
\end{equation}
which 
clearly converges to 0 as $t\rightarrow\infty$. This confirms \eqref{pozno}.

To get the first integral 
on the right of \eqref{drugadalmatinska}
we simply used that $\beta\geqslant\alpha\geqslant t^2/(2s)=:v$. And for the second one, we estimated 
$
\beta\alpha^{\frac n2-1}\geqslant \beta^{\frac12}\alpha^{\frac{n-1}2}\geqslant\sqrt{e^vv^{n-1}/2}\,,
$
and used that $v^3e^{-v}\leqslant C$.

\bigskip
We have not yet finished the estimates of $\rim{I}$ and $\rim{II}$. We still need to consider the terms in \eqref{marcia} corresponding to partial derivatives of $Q$ with respect to $\zeta$ and $\eta$.

Let $v=v(x,t)$ be as in \eqref{detinec}. Recall that the partial derivatives of $Q$ were estimated in \eqref{brcko}.  
Therefore, to estimate $\int_{\R^n} |\frac{\pd Q}{\pd\zeta}(v)||t\frac{\pd P_t f}{\pd t}(x)|\,dx$ and
 $\int_{\R^n} |\frac{\pd Q}{\pd\eta}(v)||t\frac{\pd P_t g}{\pd t}(x)|\,dx$ we need to estimate
\begin{equation*}
\label{A}
A:=\int_{\R^n} \max((P_t |f|)^{p-1}, P_t |g|)\,\bigg|t\frac{\pd P_t f}{\pd t}\bigg|\,dx
\end{equation*}
and
\begin{equation*}
\label{B}
B:=\int_{\R^n} (P_t |g|)^{q-1}\bigg|t\frac{\pd P_t g}{\pd t}\bigg|\,dx\,.
\end{equation*}

Let us prove first that
 \begin{equation}
 \label{infty}
 \lim_{t\rightarrow \infty} A = 0,\hskip 20pt \lim_{t\rightarrow \infty} B = 0\,.
 \end{equation}

To do that recall the estimate \eqref{lance} of the Hermite Poisson kernel.
It implies,
for a function $\f\in L^1$, 
\begin{equation*}
\label{ptf}P_t|\f| \,(x) \leqslant \frac{C_1}{t^n}\nor{\f}_{L^1(\R^n)} 
\end{equation*}
uniformly in $x\in\R^n$.

Now we are ready to prove \eqref{infty}. The previous inequality together with \eqref{brickroad} implies
$$
\aligned
A&= \int_{\R^n}\max((P_t|f|)^{p-1},P_t|g|)(x)\cdot\bigg|t\frac{\pd P_tf}{\pd t} (x)\bigg|\,dx\\
& \leqslant
\frac{C(n)}{t^n}(\|f\|_{1}^{p-1}+\|g\|_{1})\int_{\R^n}\int_{\R^n}\frac{t}{(|x-y|^2 + t^2)^{\frac{n+1}{2}}} |f(y)|\,dy\,dx\\
& = \frac{C'(n)}{t^n}(\|g\|_{1} + \|f\|_{1}^{p-1})\|f\|_{1} \,.
\endaligned
$$ 
This obviously tends to 0 as $t\rightarrow \infty$. 
The same with $\lim_{t\rightarrow \infty} B$. So \eqref{infty} is proved which means that $\rim{I}$ disappears.

\vspace{.1in}

We are left with the task of proving
$
 \lim_{t\rightarrow 0} A = 0\text{ and }\lim_{t\rightarrow 0} B = 0\,.
$
Let us estimate $\lim_{t\rightarrow 0} B$,  for example. From \eqref{lance} it follows that $\|P_t|g|\|_{\infty}\leqslant C(n)\|g\|_{\infty}$.
Therefore,
$$
B 
\leqslant
C(n)\|g\|_{\infty}^{q-1}
\int_{\R^n} \bigg|t\frac{\pd P_t g}{\pd t}(x)\bigg|\,dx\,.
$$
In order to estimate the integral on the right, the formula \eqref{kura} is not enough. Instead, let us denote 
$$
\Phi_n(x,y,t)=t\frac{\pd P_t(x,y)}{\pd t}  \,.
$$
This is the integral kernel of the operator 
\begin{equation}
\label{intker}
\Lambda_t:\f\mapsto t\frac{\pd P_t \f}{\pd t}\,. 
\end{equation}
From \eqref{armstrong} and \eqref{zapiski} we get 
\begin{equation}
\label{tinafey}
\Phi_n(x,y,t)=C(n)\int_0^\infty\psi'_{x,y}(u)\,\frac t{\sqrt{u}}\,e^{-\frac{t^2}{4u}}\,du\,,
\end{equation}
where $\psi_{x,y}(u)=K_u(x,y)$, as on page \pageref{sejo}.
Now
$$
\int_{\R^n} \bigg|t\frac{\pd P_t g}{\pd t}(x)\bigg|\,dx\leqslant 
\int_{\R^n}\int_{\R^n}|\Phi_n(x,y,t)||\f(y)|\,dy\,dx\,.
$$
From \eqref{tinafey} we see that the integrand on the right goes to zero pointwise as $t\rightarrow 0$. 
On the other hand, we can apply \eqref{heatderiv} to estimate $|\psi'_{x,y}(u)|$ and, consequently, \eqref{tinafey}.
The calculation which unfolds shows that, for all $t>0$, the function $(x,y)\rightarrow|\Phi_n(x,y,t)||\f(y)|$ has a majorant from $L^1(\R^n\times\R^n)$ which is independent on $t$.
This means we are entitled to use the dominated convergence theorem which gives
\begin{equation}
\label{kisa}
\lim_{t\rightarrow 0}\int_{\R^n} \bigg|t\frac{\pd P_t g}{\pd t}(x)\bigg|\,dx=0
\end{equation}
and so $\lim_{t\rightarrow 0}B=0$.
The same with $A$. 
We finally proved that $\rim{II}\leqslant 0$. 
\qed

\bigskip
\noindent{\bf Remark.}
Recall that $\Lambda_t$ were defined in \eqref{intker}. Thus \eqref{kisa} can be reformulated as $\lim_{t\rightarrow\infty}\nor{\Lambda_tg}_1=0$. 
There is also an alternative way to {\sl prove} \eqref{kisa}. First notice the statement is trivial if $g$ belongs to ${\cal V}$, the space of linear combinations of Hermite functions 
(in order to see the definition of the latter the reader is prompted to move to page \pageref{mika}). This follows from observing that $\Lambda_t=tP_t\sqrt{L}$ and applying \eqref{cluytens}, since it implies that $\sqrt L$ preserves ${\cal V}$. 
Now take arbitrary compactly supported $g$ and $\e>0$. By Lemma \ref{väyrynen} we can find $\tilde g\in{\cal V}$ such that $\nor{g-\tilde g}_1<\e$. 
Since, by \eqref{brickroad}, the kernels of $\Lambda_t$ are majorized, up to some constant $C$, by the usual Poisson kernel, we have that $\nor{\Lambda_t}_{B(L^1)}\leqslant C$ uniformly in $t$. Thus $\nor{\Lambda_t(g-\tilde g)}_1<C\e$. As $\tilde g\in{\cal V}$, there is $\delta >0$ such that for $0<t<\delta$ we have $\nor{\Lambda_t\tilde g}_1<\e$.
Therefore for such $t$ we conclude $\nor{\Lambda_t  g}_1<(C+1)\e$.

\subsection{Proof of Theorem \ref{bilinher}}

Basically we are done already. Note that, for $p\geqslant 2$, Lemma \ref{desetikvartet} and Proposition \ref{manco} together give
$$
\int_0^\infty\int_{\R^n}\nor{P_t f(x)}_*\nor{ P_t g(x)}_*\,dx\,t\,dt\leqslant C(p-1)^{}(\nor{f}_p^p+\nor{g}_q^q)
$$
with $C=48$. The same inequality obviously also holds for $1<p\leqslant 2$ if instead of $p-1$ we take $q-1$.
Now replace $f$ by $\lambda f$ and $g$ by $\lambda^{-1}g$ whereupon take the minimum in $\lambda>0$. While the left-hand side does not change, we get $C(p-1)^{}p^{1/p}q^{1/q}\nor{f}_p\nor{g}_q$ on the right-hand side. Since $p^{1/p}q^{1/q}\leqslant 2$, we obtain the desired statement of Theorem \ref{bilinher}. 
\qed

\subsection{Schr\"odinger operators with positive potential}

Following the steps in the proof of Theorem \ref{bilinher} we readily prove the next result. 

\begin{thm}
Let $V$ be a non-negative function on $\R^n$ satisfying properties \eqref{heatkato}--\eqref{emb} from page \pageref{heatkato}. 
There is an absolute constant $C>0$ such that for 
arbitrary natural numbers $M,N,n$,  
any pair $f:\R^n\rightarrow\C^M$ and $g:\R^n\rightarrow\C^N$ of $C_c^\infty$ test functions and any
$p>1$ we have
$$
\int_0^\infty\int_{\R^n}\nor{P_t f(x)}_*\nor{ P_t g(x)}_*\,dx\,t\,dt\leqslant C(p^*-1)\nor{f}_p\nor{g}_q\,.
$$
The constant $C$ 
only depends on the constant $C_0$ from the property \eqref{emb}.
\end{thm}

Let us outline the reasons why the properties \eqref{heatkato}--\eqref{emb} were needed. First, we want to make sure that the function $v(x,t)$ really maps into the domain of our Bellman function. For that reason we need that $|P_t\f(x)|^p\leqslant P_t|\f|^p(x)$. This happens if the Poisson kernel of $L$ defines a sub-probability measure at any level. But, owing to the subordination formula \eqref{assai}, it is enough to have that for the {\sl heat} kernel associated to $L$. This is exactly property \eqref{heatkato}.
Properties \eqref{gradest} and \eqref{lim}  are used to justify the estimate of the integral $-\int L'b$ from above, see sections \ref{ibp} and \ref{pritvoriseplaninata}. 
Property \eqref{emb} replaces section \ref{x^2}.

As for the lower estimates, it was already noted in a remark on page \pageref{busno1} that they are completely independent of the choice of potential $V$.

\section{Riesz transforms}

In this section we apply our embedding theorem to obtain estimates of Riesz transforms associated to Hermite operator.
Let us first introduce the necessary objects.


{\it Hermite functions} $h_m$, $m\in\N_0=\N\cup\{0\}$, are on $\R$ defined as
$$
h_m(x)=c_m(-1)^me^{x^2\over 2}{d^m\over dx^m}e^{-x^2}\,
$$
taken with the $L^2(\R)$ normalization
$
c_m={\displaystyle(2^mm!\sqrt{\pi})^{-1/2}}\,.
$

If $\alpha=(\alpha_1,\hdots,\alpha_n)\in\N_0^n$, 
then the $Hermite$ $function$ \label{mika} on $\R^n$ is given by
$$
h_\alpha:=h_{\alpha_1}\otimes\hdots\otimes h_{\alpha_n}\,.
$$
Next we provide some argumentation as to why it is convenient to take 
linear combinations of 
$h_\alpha$'s
as the family of test functions.

\begin{lema}
\label{väyrynen}
The space {\rm Lin}\,$\mn{h_\alpha}{\alpha\in\N_0^n}$ is dense in $L^p(\R^n)$ for $1\leqslant p<\infty$.
\end{lema}
\dok
Take $f\in C_c^\infty(\R^n)$ and define coefficients $\hat{f}(\alpha)$ as
$$
\hat{f}(\alpha)=\int_{\R^n}f(x)h_\alpha(x)\,dx\,.
$$
Consider
$$
S_Nf:=\sum_{|\alpha|\leqslant N}\hat{f}(\alpha)h_\alpha\,.
$$
We can repeat the proof of Lemma 5.4.1 in \cite{T} to show that $S_Nf$ converge to $f$ in the $L^p$ norm.
\qed

\bigskip
\noindent
{\bf Remark.}
The previous sentence is not true for arbitrary $f\in L^p$. In fact, a well-known theorem by Askey and Wainger states that already when $n=1$, this {\it is} the case if and only if $p^*<4$.

\bigskip
Recall that $\z_j$ and $\z_j^*$ were introduced in \eqref{mole}.
By \cite[1.1.30]{T},
\begin{equation}
\label{ludovic}
\aligned
\z_jh_\alpha & =\sqrt{2(\alpha_j+1)}\,h_{\alpha_1}\otimes\hdots\otimes h_{\alpha_{j}+1}\otimes\hdots\otimes h_{\alpha_n}\\
& =\sqrt{2(\alpha_j+1)}\,h_{\alpha+e_j}
\endaligned
\end{equation}
and similarly
\begin{equation*}
\label{vaillant}
\z_j^*h_\alpha=\sqrt{2\alpha_j}\,h_{\alpha-e_j}\,.
\end{equation*}
Together with \eqref{nikolaeva} this implies
\begin{equation}
\label{cluytens}
L h_\alpha=(2|\alpha|+n)h_\alpha\,,
\end{equation}
where $|\alpha|=\alpha_1+\hdots+\alpha_n$.

\subsection{Spectral multipliers}
\label{sm}

Another tool we need in order to treat the Riesz transforms are spectral multipliers. They are defined as follows. 

Let $\Psi$ be a bounded complex function on $\N$. In view of \eqref{cluytens} it is natural to define 
\begin{equation*}
\Psi(L):=\sum_{m=0}^\infty\Psi(2m+n)\,{\cal P}_m\,,
\end{equation*}
where ${\cal P}_m$ is the projector onto the subspace of $L^2(\R^n)$ generated by all Hermite functions $h_\alpha$ with $|\alpha|=m$. 

We are interested in $L^p$ boundedness of such operators. The sheer boundedness of $\Psi$ does not guarantee that (unless $p=2$). 
Certain sufficient conditions are given by the 
multiplier theorems for Hermite expansions 
due to Mauceri \cite{M} and Thangavelu \cite[Theorem 4.2.1]{T}. They imposed Marcinkiewicz-H\"ormander-Mihlin-type conditions on their multipliers. 
Results in the same spirit were also obtained for closely related operators such as Weyl multipliers \cite{M}, \cite{T2} and multipliers associated with the twisted Laplacian \cite{N}.
We should also mention the paper \cite{DOS} which itself contains many further references.
But for our purpose we need more -- we want estimates independent of $n$ and $p$. 
They are provided by 
Theorem \ref{multiplier} which is proven below. It basically confirms the assertion made in \cite{DV}, see Remark 3.2 there, namely, that the method exposed in \cite{DV} is only dependent on successful treatment of the corresponding spectral multipliers. The latter, in turn, should follow from ``non-singularity" of the spectrum of the underlying differential operator (in our case, $L$).


\bigskip
\noindent
{\bf Proof of Theorem \ref{multiplier}.}
Let $\Phi(z)=\Psi(1/z)$. The assumption on $\Psi$ can be restated to say that $\Phi$ is analytic in a neighbourhood of 0. If $R$ is the radius of convergence of its power series expansion around 0, let $\rho=1/R$. Thus 
$\Psi$ is analytic in $\{|z|>\rho\}$. 

For the sake of convenience assume that $\Phi(0)=0$; it is trivial to remove this restriction. Therefore $\Phi$ can be expanded as $\Phi(w)=\sum_{j=1}^{\infty}c_jw^j$, provided that $|w|<R$. 
The Cauchy-Hadamard formula gives
\begin{equation}
\label{parsifal}
\rho=\limsup_{j\rightarrow\infty}|c_j|^{1/j}\,. 
\end{equation}

Suppose first that $n>\rho$. In that case we can write 
$
\Psi(L)=\Phi(L^{-1})=\sum_{j=1}^{\infty}c_jL^{-j}\,.
$ 
Choose 
$a>\rho$
which is also strictly smaller than any integer that exceeds $\rho$. In other words, fix $a\in(\rho,[\rho]+1)$.
We have the formula 
\begin{equation*}
\label{ergen-deda}
L^{-{j}} = \frac{1}{(j-1)!}\int_0^{\infty} t^{j-1} e^{-tL}\, dt\,.
\end{equation*}
Then
\begin{equation*}
\|L^{-{j}}\|_{p\rightarrow p} \leqslant \frac{1}{(j-1)!}\int_0^{\infty} t^{j-1} \|e^{-tL}\|_{p\rightarrow p}\, dt\,.
\end{equation*}
Since $a>\rho$, it follows from \eqref{parsifal} that there is $C=C(a)>0$ such that $|c_j|\leqslant C a^j$ for all $j\in\N$.
Consequently,
\begin{equation}
\label{pilence}
\nor{\Phi(L^{-1})}_{p\rightarrow p}\leqslant 
Ca\int_0^{\infty} e^{at} \|e^{-tL}\|_{p\rightarrow p}\, dt\,.
\end{equation}
Hence in order to proceed we must estimate the norm of $e^{-tL}$ on $L^p$.

Heat semigroup $e^{-tL}$ is known to be contractive in every $L^p$.
Moreover, in $L^2$ it is very contractive, in the sense that
\begin{equation}
\label{stanke}
 \|e^{-tL}\|_{2\rightarrow 2} = 
 e^{-nt}\,,
\end{equation}
just because the smallest eigenvalue of $L$ is $n$ if we are in a $n$-dimensional space.
 
On the other hand, we have by \eqref{8-2} that
\begin{equation}
\label{kalimankov}
\nor{e^{-tL}}_{\infty\rightarrow\infty}
\leqslant 
{1\over (2\pi)^{n/2}}\sup_{x\in\R^n}\int_{\R^n}K_t(x,y)\,dy\,,
\end{equation}
where $K_t$ is as in \eqref{harlech}. By taking $\f\equiv 1$ in \eqref{8-2} we see that we actually have an equality in \eqref{kalimankov}.
A calculation, very similar to the one from the first part of Section \ref{pritvoriseplaninata}, shows that 
\begin{equation}
\label{abc}
{1\over (2\pi)^{n/2}}\int_{\R^n}K_t(x,y)\,dy=(\cosh2t)^{-n/2}\,e^{-\frac{|x|^2}{2\coth2t}}\,,
\end{equation}
hence we can continue \eqref{kalimankov} as 
\begin{equation}
\label{denkov}
\nor{e^{-tL}}_{\infty\rightarrow\infty}
=(\cosh2t)^{-n/2}
\,.
\end{equation}
Note that exactly the same identity 
is valid on $L^1$.

Complex interpolation between \eqref{stanke} and \eqref{denkov} yields, for arbitrary $p\in[1,\infty]$,
\begin{equation*}
\nor{e^{-tL}}_{p\rightarrow p}
\leqslant e^{-2nt/p^*}(\cosh 2t)^{-\frac{n}2\gamma(p)}\,,
\end{equation*}
where 
$
\gamma(p)=1-2/{p^*}\,.
$
By applying it in \eqref{pilence} we get
\begin{equation}
\label{todor}
\aligned
\nor{\Phi(L^{-1})}_{p\rightarrow p} 
& \leqslant 
Ca\int_0^{\infty} e^{(a-n)t}[e^{-2t}\cosh 2t]^{-\frac{n}2\gamma(p)}\, dt\,.
\endaligned
\end{equation}
This integral converges if and only if $n>a$, that is, if and only if $n>\rho$. 
In order to estimate it from above, we first estimate the $\cosh$ part of the integrand as follows:
\begin{equation}
\label{zavidovici}
\cosh 2t\geqslant\max\{1,\,e^{2t}/2\}=
\left\{
\begin{array}{lcr}
1 & ; & t\leqslant t_0\\
e^{2t}/2 & ; & t\geqslant t_0
\end{array}
\right.\,,
\end{equation}
where 
$
t_0=\log\sqrt2\,. 
$
Next we write the integral in \eqref{todor} as
$$
\int_0^{\infty}=\int_0^{t_0}+\int_{t_0}^{\infty}
$$
and in each of the integrals on the right apply the appropriate estimate from \eqref{zavidovici}. The resulting integrals can be explicitly calculated. One obtains
\begin{equation}
\label{zusatz}
\frac1{Ca}\nor{\Phi(L^{-1})}_{p\rightarrow p}\leqslant \frac{\sqrt{2}^\alpha-1}\alpha+\frac{\sqrt{2}^\alpha}\beta\,,
\end{equation}
where $\alpha=a-2n/p^*$ and $\beta=n-a$. This estimate is valid as long as $\beta>0$ (in the case of $\alpha=0$ one has the limiting expression, i.e. $\log\sqrt2+1/\beta$). 
The right-hand side of \eqref{zusatz} is increasing in $\alpha$. But 
obviously $\alpha<a$, which means that the expression in \eqref{zusatz} is uniformly bounded as $p$ ranges over $[1,\infty]$ and $n$ ranges over integers bigger than $\rho$. 
The assumption $n>\rho$ was used here for the second time, namely, to estimate $\beta\geqslant [\rho]+1-a>0$.

\medskip
Now let us consider the case when $n\leqslant\rho$. Assume first that $n=1$. We are indebted to Adam Sikora for showing us how to proceed in this case. 

We would like to show that the operator $(L+I)^{-1/2}$ maps boundedly $L^2\rightarrow L^\infty$ (then, by duality, we also get $L^1\rightarrow L^2$ boundedness). Indeed,
for a Schwarz function $u$ we have
\begin{equation}
\label{bogen}
\aligned
|u(x)|&\leqslant\int_\R|\widehat u|=\int_\R|\widehat u(x)|(1+|x|^2)^{1/2}(1+|x|^2)^{-1/2}\,dx\\
&\leqslant C \Big(\int_\R|\widehat u(x)|^2(1+|x|^2)\,dx\Big)^{1/2}=C\sk{(|x|^2+1)\widehat u}{\widehat u}^{1/2}\\
&=C\sk{[(-\Delta+I)u]\widehat\ }{\widehat u}^{1/2}=C\sk{(-\Delta+I)u}{u}^{1/2}\leqslant C\sk{(L+I)u}{u}^{1/2}\,.
\endaligned
\end{equation}

Conversely, let us show that $(L+I)^{-1/2}$ is also bounded as an operator $L^\infty\rightarrow L^2$. Since this is a self-adjoint operator, it suffices to show that it is bounded $L^2\rightarrow L^1$. 
But this simply follows from the H\"older's inequality:
\begin{equation}
\label{kreuzer}
\nor{u}_1^2\leqslant C \nor{(|x|^2+I)^{1/2}u}_2^2=C\sk{(|x|^2+I)u}{u}\leqslant
C\sk{(L+I)u}{u} \,.
\end{equation}

At the beginning of the proof we assumed that $\Phi(z)=0$. This quickly implies that 
\begin{equation}
\label{hajro}
|\Psi(\lambda)|\leqslant \frac{C}{\lambda +1}\,
\end{equation}
for some $C>0$.
Therefore 
\begin{equation*}
\aligned
\nor{\Psi(L)}&_{\infty\rightarrow \infty}
\leqslant
\nor{(L+I)^{-1/2}}_{2\rightarrow \infty}
\nor{\Psi(L)(L+I)}_{2\rightarrow 2}
\nor{(L+I)^{-1/2}}_{\infty\rightarrow 2}
\,.
\endaligned
\end{equation*}
Thus $\Psi(L)$ is bounded $L^\infty\rightarrow L^\infty$. Now the Riesz-Thorin theorem implies that there is an absolute $C>0$, such that $\nor{\Psi(L)}_{p\rightarrow p}\leqslant C$ for all $p\geqslant 2$. Since $\Psi(L)^*=\overline\Psi(L)$, and since $\overline\Psi$ also satisfies \eqref{hajro}, we get $\nor{\Psi(L)}_{p\rightarrow p}\leqslant C$ for 
$1< p\leqslant 2$, as well, while boundedness on $L^1$ can be proven directly as above. 

Basically the same proof is valid for arbitrary $n\in\N$. The difference is that in general one needs to take the $n/2$-th power of $x^2+1$ in order to run the H\"older's inequality in \eqref{kreuzer}. 
So one should deal with $(L+I)^n$. 
But when $n>1$ one cannot, as in \eqref{bogen} and \eqref{kreuzer}, simply discard the terms $|x|^2$ or $-\Delta$ in the lower estimate of $\sk{(L+I)^nu}{u}$. However, the $L^2\rightarrow L^1$ boundedness of $S:=(L+I)^{-n/2}$ is for arbitrary $n$ provided by the estimate (7.11) from \cite{DOS}. To obtain the $L^2\rightarrow L^\infty$ boundedness we use the Fourier transform. It convenes us to define it 
on $\R^n$ as 
$$
\hat{f}(\xi)=(2\pi)^{-n/2}\int_{\R^n}f(x)\,e^{-i\sk{x}{\xi}}dx\,.
$$
For then $\widehat{-\Delta f}(\xi)=|\xi|^2\hat f(\xi)$ and $\widehat{|x|^2f}(\xi)=-\Delta\hat f(\xi)$, thus $\widehat{Lu}=L\widehat u$, i.e. the Fourier transform commutes with $L$.
Consequently,
$$
\nor{Su}_\infty\leqslant\nor{\widehat{Su}}_1=\nor{S\hat{u}}_1\leqslant C\nor{\widehat u}_2=C\nor{u}_2\,.
$$

Finally, these modifications call for a suitably sharper estimate in \eqref{hajro}, namely $|\Psi(\lambda)|\leqslant C(\lambda +1)^{-n}$, in order to repeat the calculation which follows it. But this can be easily achieved, since we may assume without loss of generality that the first $n$ derivatives of $\Phi$ at zero vanish.  

All this shows that $\nor{\Psi(L)}_{B(L^p(\R^n))}\leqslant C(n)$. Now just take maximum of the constant, which appeared in the estimates for $n>\rho$, and all $C(n)$ for $1\leqslant n\leqslant \rho$. This is our absolute constant. 
\qed

\subsection{Proof of Corollaries \ref{RieszH} and \ref{schwertleite}}

If viewed correctly, these 
are consequences of Theorems \ref{bilinher} and \ref{multiplier}. 
The connection between them 
will be established through the following two formulas:
\begin{equation}
\label{f1}
\sk{R_j f}{g} =\int_{0}^{\infty} \Sk{\z_j P_t {\cal O} f}{\frac{\pd}{\pd t}P_t g}_{L^2(\R^n)}\, t\, dt\,
\end{equation}
and
\begin{equation}
\label{f2}
\langle R_j^*f, g\rangle =\int_{0}^{\infty}\Sk{\z_j^* P_t {\cal O}^* f}{\frac{\pd}{\pd t}P_t g}_{L^2(\R^n)}\, t\, dt\,.
\end{equation}
Here ${\cal O}$ and ${\cal O}^*$ are operators in $L^p(\R^n)$, $1<p<\infty$, hopefully bounded independently of dimension $n$.
In order to calculate these operators we test the formulas on Hermite functions. More precisely, take
$$
f:=L^\frac12h_\alpha\,.
$$
Then \eqref{f1} becomes
\begin{equation}
\label{jagode}
\sk{\z_j h_\alpha}{g} =-\int_{0}^{\infty} \Sk{ L^{\frac12}P_t\z_j P_t {\cal O}  L^{\frac12}h_\alpha}{g}\, t\, dt\,.
\end{equation}

Write formally
\begin{equation}
\label{O!}
{\cal O}h_\alpha=\sum_{\beta\in\N_0^n}o_{\alpha\beta}h_\beta\,.
\end{equation}
By using \eqref{ludovic}, \eqref{cluytens} and \eqref{O!}, we formally calculate
$$
L^{\frac12}P_t\z_j P_t {\cal O}  L^{\frac12}h_\alpha=\sqrt{2|\alpha| + n}\sum_{\beta\in\N_0^n}o_{\alpha\beta}
 e^{-t\sqrt{2|\beta| +n}}\lambda_{|\beta|+1}(t)\sqrt{2(\beta_j+1)}\,h_{\beta+e_j}\,.
$$
Here
$$
\lambda_k(t)=\sqrt{2k+n}\,e^{-t\sqrt{2k+n}}
$$
Together with \eqref{jagode} this means that we can take $o_{\alpha\beta}=0$ if $\alpha\not =\beta$ and we have a formula for the coefficients $o_{\alpha\alpha}$, which we can denote by $o_\alpha$:
$$
{1\over o_\alpha}=-\int_{0}^{\infty}\, \lambda_{|\alpha|}(t)\lambda_{|\alpha|+1}(t)\,t\, dt\,.
$$
Thus $o_\alpha$ and $o_\alpha^*$ actually depend on $|\alpha|$ only, so if we denote $m=|\alpha|$, we may write $o_\alpha=o_m$, $o_\alpha^*=o_m^*$. We get
$$
o_m   =-\frac{(\sqrt{2m+n}+\sqrt{2m+n+2})^2}{\sqrt{2m+n}\sqrt{2m+n+2}}\,,\hskip 17pt
o_m^* =-\frac{(\sqrt{2m+n}+\sqrt{2m+n-2})^2}{\sqrt{2m+n}\sqrt{2m+n-2}}\,
$$
or, equivalently,
\begin{equation*}
\label{mustafov}
o_m=-\Psi(2m+n)\hskip 15pt \text{and}\hskip 15pt o_m^*=o_{m-1}\,,
\end{equation*}
where
\begin{equation}
\label{griffi}
\Psi(k)=\frac{(\sqrt{k}+\sqrt{k+2})^2}{\sqrt{k}\sqrt{k+2}}\,.
\end{equation}
Consequently,
\begin{equation*}
\label{bussueil}
{\cal O}=\sum_{m\in\N_0}o_m{\cal P}_m\hskip 15pt \text{and}\hskip 15pt
{\cal O}^*=\sum_{m\in\N_0}o_m^*{\cal P}_m\,,
\end{equation*}
recalling that 
${\cal P}_m$ is the projector onto the subspace of $L^2(\R^n)$ generated by all Hermite functions $h_\alpha$ with $|\alpha|=m$. 

Note that $o_0^*$ is defined if $n\geqslant 3$. So if $n=1,2$ we have to correct formula \eqref{f2} as $o_0^*:=0$ and
$$
\langle R_j^*f, g\rangle =\sk{R_j^*{\cal P}_0f}{g}+\int_{0}^{\infty}\Sk{\z_j^* P_t {\cal O}^* f}{\frac{\pd}{\pd t}P_t g}\, t\, dt\,.
$$
But \eqref{cluytens} implies that $L^{\frac12}h_0=\sqrt{n}\,h_0$, therefore
$$
R_j^*{\cal P}_0f=\z_j^*L^{-\frac12}\sk{f}{h_0}h_0=\sk{f}{h_0}{1\over\sqrt{n}}\z_j^*h_0=0\,.
$$
We actually proved that for any $n$ we can take $o_0^*=0$ and
\eqref{f2} remains valid.

\bigskip
\noindent
{\bf Remark.}
It does not come as a surprise that the formul\ae\ for $o_m,o_m^*$ are very similar to those in \cite[p. 183]{DV}. See also Remark 3.2 in the same paper.

\bigskip
By applying Theorem \ref{multiplier} we immediately get the following result.

\begin{prop}
\label{pivka}
For all $p\in[1,\infty]$ and all $n\in\N$, operators ${\cal O}$ and ${\cal O}^*$ are bounded on $L^p(\R^n)$ with constants independent of $n$ and $p$.
\end{prop}
\dok
Indeed, ${\cal O}=\Psi(L)$ with $\Psi$ as in \eqref{griffi}, while ${\cal O}^*=\Psi^*(L)$, where $\Psi^*(k)=\Psi(k-2)$ for $k>2$ and $\Psi^*(1)=\Psi^*(2)=0$.\qed

\bigskip
Let us show how \eqref{f1} and \eqref{f2} help to estimate Riesz transforms. 
Take $m\in\N$ and let $f=(f_1,\hdots,f_m)$ be a $\C^m-$valued test function on $\R^n$. By $R_jf$ we will understand $(R_jf_1,\hdots,R_jf_m)$; similarly for $R_j^*$.
Let also ${\mathcal R}f=(R_1f,\hdots,R_nf,R_1^*f,\hdots,R_n^*f)$. This is a function with values in $(\C^m)^{2n}$. 
Thus we can think of ${\mathcal R}f$ as a matrix function with entries $R_jf_k$ and $R_j^*f_k$, where $j=1,\hdots,n$ and $k=1,\hdots,m$. Therefore
$$
\nor{{\mathcal R}f}_{p}^p
=\int_{\R^n}\Big(\sum_{j,k}|R_jf_k(x)|^2+|R_j^*f_k(x)|^2\Big)^{p/2}\,dx
=\int_{\R^n}\nor{{\mathcal R}f(x)}_{HS}^p\,dx\,,
$$
where $\nor{\cdot}_{HS}$ stands for the usual Hilbert-Schmidt norm.
\begin{prop}
\label{sieglinde}
Under the above notation, 
$$
\nor{{\mathcal R}f}_{p}\leqslant C(p^*-1)\nor{f}_{p}\,.
$$
for some absolute $C>0$ and all $n\in\N$, $p>1$.
\end{prop}

Note that Corollary \ref{RieszH} is just a special case (when $m=1$) of this propositon, while Corollary \ref{schwertleite} follows by applying it repeatedly. Therefore it is only left for us to prove Proposition \ref{sieglinde}.

\bigskip
\dok
Observe that 
$$
\nor{{\mathcal R}f}_{p}=\sup|\sk{(R_1f,\hdots,R_{n}f,R_1^*f,\hdots,R_{n}^*f)}{(g_1,\hdots,g_{2n})}|\,,
$$
the supremum being taken over all $g=(g_1,\hdots,g_{2n})$ with $L^q$ norm (in the appropriate space) equal to one, where each $g_j$ is a function $\R^n\rightarrow\C^m$. 

Now by \eqref{f1} and \eqref{f2},
$$
\aligned
& |\sk{(R_1f,\hdots,R_{n}f,R_1^*f,\hdots,R_{n}^*f)}{(g_1,\hdots,g_{2n})}|=
\\
& \bigg|\int_0^\infty\int_{\R^n}
\sum_{j=1}^{n}\Big(
\sk{\z_j P_t {\cal O} f(x)}{\pd_tP_t g_j(x)}_{\C^m}+
\sk{\z_j^* P_t {\cal O}^* f(x)}{\pd_tP_t g_{n+j}(x)}_{\C^m}\Big)\,dx\,t\,dt\bigg|\,.
\endaligned
$$
Here, as before, by ${\cal O} f$ we mean $({\cal O} f_1,\hdots,{\cal O} f_m)$; similarly for $\z_j P_t$, $\z_j^* P_t{\cal O}^*$ and $\pd_tP_t$. By the Cauchy-Schwarz inequality we get
\begin{equation}
\label{kudapojdu}
\leqslant \int_0^\infty\int_{\R^n}
\bigg(\sum_{j=1}^{n}\Big|\z_j P_t {\cal O} f(x)\Big|^2+\sum_{j=1}^{n}\Big|\z_j^* P_t {\cal O}^* f(x)\Big|^2\bigg)^{1/2}
\bigg(\sum_{j=1}^{2n}\Big|\pd_tP_t g_j(x)\Big|^2\bigg)^{1/2}\,dx\,t\,dt\,.
\end{equation}
We continue by a raw estimate
\begin{equation}
\label{zaschischaya}
\leqslant \sqrt{2}\int_0^\infty\int_{\R^n}
\Big(\nor{P_t {\cal O} f(x)}_*+\nor{P_t {\cal O}^* f(x)}_*\Big)
\nor{P_t g(x)}_{*}
\,dx\,t\,dt 
\end{equation}
whereupon 
Theorem \ref{bilinher} 
yields
\begin{equation}
 \label{waltraute}
 \leqslant C(p^*-1)
 \Big(\nor{{\cal O} f}_p+\nor{{\cal O}^* f}_p\Big)
 \nor{g}_q\,.
 \end{equation}

In Proposition \ref{pivka} we proved 
$\nor{{\cal O}}_{B(L^p(\R^n))}\leqslant C$ 
with some absolute $C>0$ 
and the same for $\nor{{\cal O}^*}_{B(L^p(\R^n))}$.
From the theorem of Marcinkiewicz and Zygmund (see, for example, \cite{G}) it follows that the same bounds apply to the $l^2-$valued extensions of ${\cal O}$ and ${\cal O}^*$ that appear in \eqref{waltraute}.
Therefore the proof of Proposition \ref{sieglinde} (and consequently of Corollaries \ref{RieszH} and \ref{schwertleite}) is complete. 
\qed

\subsection{Heisenberg groups}
\label{Heisenberg}

In this section we review some of the well-known ties between the Hermite operator and the Heisenberg groups. 
We present the result by Coulhon, M\"uller and Zienkiewicz \cite{CMZ} establishing the dimension-free $L^p$ boundedness of the Heisenberg-Riesz transform. Their estimate involves the $L^p$ norms of the Hilbert transform $\cH$  along the parabola in the plane. 
In order to push the approach of \cite{CMZ} to its limit we devote special attention to 
summarizing
a series of deep results by Seeger, Tao and Wright obtained in recent years which together give (almost) sharp $L^p$ estimates of $\cH$.

For the purpose of keeping this section as self-contained as possible we need to start with the basic definitions. 
They can be found in most of the introductory texts on Heisenberg groups, for example \cite{S}, \cite{T1} or \cite{Ta}.

\bigskip
By the Heisenberg group $\He^n$ we understand $\R^{2n+1}$ endowed with the multiplication
$$
(u,z)\cdot(u',z')=\Big(u+u',z+z'+2\sum_{j=1}^n(x_j'y_j-x_jy_j')\Big)
\,.
$$
Here $u=(x_1,y_1,\hdots,x_n,y_n)\in\R^{2n}$, $z\in\R$, and similarly for $(u',z')$.

Let $\lambda:\He^n\rightarrow {\rm Lin}({\cal S}(\He^n))$ be the left regular representation of $\He^n$, i.e.  $\lambda: h\mapsto \lambda_h$ where $\lambda_hf(v)=f(h^{-1}v)$. An operator $S\in {\rm Lin}({\cal S}(\He^n))$ is said to be {\it left-invariant} if $S\circ \lambda_h=\lambda_h\circ S$ for every $h\in\He^n$.

The underlying Lie algebra $\h^n$ of all left-invariant vector fields on $\He^n$ is generated by the fields 
$$
\aligned
\label{etot}
X_{2k-1} & = & \frac{\partial}{\pd x_k}+2y_k\,\frac{\partial}{\pd z}\\
X_{2k} & = & \frac{\partial}{\pd y_k}-2x_k\,\frac{\partial}{\pd z}
\endaligned
$$
for $k\in\{1,\hdots,n\}$ and
$$
Z=\frac{\partial}{\pd z}\,.
$$ 
The Lie bracket in $\h^n$ is defined as $[U,V]\f=U(V\f)-V(U\f)$.

In \cite{CMZ} the authors deal with the Riesz transforms on the Heisenberg group, defined as 
$X_j\LL^{-1/2}$, $j=1,\hdots,2n$, where 
$$
\LL=-\sum_{j=1}^{2n}X_j^2
$$
is the sub-Laplacian on $\He^n$. They prove, with some constant $C(p)$ independent of $n$, that 
$$
\Nor{{\bf R}_{\He^n}f}_p\leqslant C(p)\nor{f}_p\,,
$$
where 
$$
{\bf R}_{\He^n} f=\Big(\sum_j|X_j\LL^{-1/2}f|^2\Big)^{1/2}\,.
$$
A feature of their proof is the use of Gaveau-Hulanicki formula for heat kernels on the Heisenberg group which is utilized for representation of the Riesz transforms (see the formula on p. 378 of \cite{CMZ} which, due to a misprint, requires a minor modification; compare with p. 374 there):
$$
(X_k\LL^{-1/2})f(u)=\int_{\He^n}\Phi_k(v)\,\tilde\cH_vf(u)\,dv\,,
$$
where $u=(x_1,y_1,\hdots,x_n,y_n,s)\in\He^n$,
$$
\Phi_k=\left\{
\begin{array}{rcl}
x_jF-i\,y_jG & ; & k=2j-1\\
y_jF+i\,x_jG & ; & k=2j
\end{array}
\right.
$$ 
and, for $v=(w,z)\in\He^n$,
$$
\aligned
F(v)=\int_{-\infty}^{\infty}{\bf p}(v,\lambda)\cosh\lambda\,d\lambda\\
G(v)=\int_{-\infty}^{\infty}{\bf p}(v,\lambda)\sinh\lambda\,d\lambda
\endaligned
$$
with 
$$
{\bf p}(v,\lambda)=-\frac1{4(2\pi)^{n+1}}\,e^{-\frac\lambda2(|w|^2\coth\lambda-i\,z)}
\Big(\frac\lambda{\sinh\lambda}\Big)^{n+1}\,.
$$
Furthermore, $\tilde\cH_v$ is the Hilbert transform along the parabola $t\mapsto \delta_t(v)$ in $\He^n$:
$$
(\tilde\cH_v f)(u)={\rm p.v.}\int_{-\infty}^{\infty} f(u\cdot\delta_t(v)^{-1})\,\frac{dt}{t}\,.
$$  
Here $\delta_t$ denotes the Heisenberg dilations $\delta_t(w,z)=(tw,t^2z)$. 
The authors observe that the norms of $\tilde\cH_v$ can be reduced to estimating the Hilbert transform along the standard parabola in $\R^2$, i.e. 
$$
(\cH f)(x,y)={\rm p.v.}\int_{-\infty}^{\infty}f(x-t,y-t^2)\,\frac{dt}{t}\,,\hskip 15pt (x,y)\in\R^2\,.
$$
By proceeding as in Section III of \cite{CMZ} one arrives at the estimate 
\begin{equation}
\label{daviti}
\nor{{\bf R}_{\He^n}f}_p\leqslant C_{q,n} \nor{\cH}_p\nor{f}_p 
\end{equation}
with
$$
C_{q,n}=[\sigma(S^n)]^{1/p} \Nor{\int_0^\infty X_1p_1(\delta_r(\omega))\,r^{2n+1}dr}_{L^q(d\sigma(\omega))}\,.
$$
Here $S^n$ stands for the unit sphere in the Kor\'anyi norm, given 
by $\nnor{(u,z)}^4=|u|^4+z^2$, while $\sigma$ is the induced surface measure on $S^n$. The Kor\'anyi norm is homogeneous with respect to the Heisenberg dilations, 
meaning that $\nnor{\delta_r\omega}=r\nnor{\omega}$ for any $\omega\in\He^n$. Furthermore, $p_1$ is the heat kernel on $\He^n$ calculated at the level 1; it is explicitly given by the formula due to Gaveau \cite{G1} and Hulanicki \cite{H}. 

Strictly speaking $C_{q,n}$ depends on $n$, but the authors devote their Section IV to showing that it 
has a majorant which 
does not. Actually, a careful examination of their proof reveals that 
it
can even be estimated from above by an {\it absolute} constant, i.e. one indepentent of both $n$ and $q$. In order to prove that $C_{q,n}\leqslant C$ for some $C>0$ and all $q>1$, $n\in\N$, let us attempt
rewriting 
the bottom line on p. 375 from \cite{CMZ} with the use of \cite[Lemma 4]{CMZ} 
but with absolute $C$ instead of $C_q$:
$$
\aligned
2\frac{\pi^{n-1/2}}{\Gamma(n-1/2)}n^{-1-q/2}\leqslant 
C^q  
& \,\pi^{nq}2^{-nq}\frac{\Gamma(n)^q}{\Gamma(n+3/2)^q\Gamma(n/2)^{2q}}
\cdot\\
\cdot
& (4\pi^{n+1/2})^{1-q}\frac{\Gamma(n/2)^{1-q}}{\Gamma(n)^{1-q}\Gamma((n+1)/2)^{1-q}}
\,.\endaligned
$$
By means of the duplication formula for the $\Gamma$ function \cite[A-5]{G} proving the above inequality is the same as proving 
$$
\frac{2^{-nq+4q+n-2}
\pi^{3(q-1)/2}}
{n^{1+q/2}}
\cdot
\frac{\Gamma(n+3/2)^q\Gamma(n/2)^{q}}
{\Gamma(n-1/2)\Gamma(n/2+1/2)^{3q-2}}
\leqslant 
C^q 
\,.
$$
The identity $\Gamma(x+1)=x\,\Gamma(x)$ 
translates this into
$$
\frac{2^{-nq+4q+n-2}
\pi^{3(q-1)/2}}
{n^{1+q/2}}
\cdot
\frac{(n+1/2)^q(n-1/2)\Gamma(n+1/2)^{q-1}\Gamma(n/2)^{q}}
{\Gamma(n/2+1/2)^{3q-2}}
\leqslant 
C^q 
\,.
$$
Due to the Stirling formula this is equivalent to proving
$$
\frac{2^{4q-3/2}
e^{1/2}
\pi^{q-1}(n+1/2)^q
(n-1/2)^{n(q-1)+1}
(n-2)^{(n-1)q/2}}
{n^{1+q/2}(n-1)^{(3q-2)n/2}}
\leqslant 
C^q 
\,,
$$
that is,
$$
2^{4q-3/2}
e^{1/2}
\pi^{q-1}
\cdot
\frac{n-1/2}n
\cdot
\left[
\frac{(n+1/2)^2}{n(n-2)}
\right]^{q/2}
\cdot
\left(
\frac{n-1/2}{n-1}
\right)^{n(q-1)}
\cdot
\left(
\frac{n-2}
{n-1}
\right)^{nq/2}
\leqslant 
C^q 
\,.
$$
From this type of expression it is clear that indeed there is an {\it absolute} constant $C>0$ 
such that for all $n\in\N$ and all $q$, the left-hand side is majorized by 
$C^q$.

All said implies the following improvement of \eqref{daviti}, i.e. the explicite estimate in the major result of \cite{CMZ} with an absolute constant $C$:
\begin{equation}
\label{podpole}
\nor{{\bf R}_{\He^n}f}_p\leqslant C\nor{\cH}_p\nor{f}_p 
\,,
\end{equation}
which is to say that the $L^p$ norm of the vector Heisenberg-Riesz transform is controlled by the $L^p$ norm of $\cH$ alone.

\subsubsection*{Orlicz spaces and the estimate of $\cH$}

The $L^p$ boundedness of $\cH$ has been studied for many years. In 1966 it was proven by Fabes \cite{F} that $\cH$ is bounded on $L^2$. Later on Stein and Wainger \cite{SW} extended this result to $L^p$ for $1<p<\infty$. 
As for $p=1$, 
the problem of determining optimal boundedness of $\cH$ has been a major challenge in the area. 
There is, namely, a long-standing  conjecture 
that $\cH$ is of weak type $(1,1)$. This question is still open. A close result in this direction, due to Seeger, Tao and Wright \cite{STW}, says that $\cH$ maps from $L\log\log L$ to weak $L^1$ space. 
Therefore, for every $\e>0$ it maps from $L\log^\e L$ to weak $L^1$. By invoking an interpolation argument as done by Tao and Wright \cite{TW} we see that it maps $L\log^{1+\e} L$ to (strong) $L^1$. 
This chain of implications is completed by a theorem of Tao \cite{Tao} which implies that $\cH$ is bounded on $L^p$ with the estimate $\nor{\cH}_p\leqslant C_\e p^{1+\e}$, where $C_\e>0$.

It is worthwhile 
putting down the amalgam of the above paragraph and
\eqref{podpole}:\\
 
{\it For every $\e>0$ there is $C_\e>0$ such that}
\begin{equation}
\label{msu_notredame}
\nor{{\bf R}_{\He^n}f}_p\leqslant C_\e p^{1+\e}\nor{f}_p  \,.
\end{equation}

Of the vast literature existing on general types of Radon transforms and their estimates we single out the work by Christ, Nagel, Stein and Wainger \cite{CNSW}.

\subsubsection*{Schr\"odinger representations}

It is left to explain why \eqref{msu_notredame} implies the same estimate for the Hermite-Riesz transforms ${\bf R}$.

This is done by means of the {\it Schr\"odinger representations} of $\He^n$.  For a fixed $\lambda\in\R\backslash\{0\}$ define the operator $\pi_\lambda: \He^n\rightarrow {\cal U}(L^2(\R^n))$ by 
$$
[\pi_\lambda(x,y,z)f](v)=e^{i\lambda(x\cdot v+\frac12 x\cdot y+\frac14 z)}f(v+y)\,.
$$
Here $x,y,v\in\R^n$, $z\in\R$ and the dot denotes the usual Euclidean scalar product in $\R^n$. This gives rise to a derived representation, let us call it $d\pi_\lambda$, of the Lie algebra $\h^n$. It is defined as follows.

Let $\gamma:[0,1]\rightarrow \He^n$ be a $C^1$ curve such that $\gamma(0)=0$ and assume 
$\Xi\in\h^n$ is given by 
$$
\Xi f(a)=\frac{d}{ds}f(a\cdot \gamma(s))\Big|_{s=0}
$$
where  $a\in\He^n$ and $f:\He^n\rightarrow \C$. 
By choosing the coordinate curves 
$$
\gamma_j(s)=(0,\hdots,0,s,0,\hdots,0)
$$ 
with $s$ on the $j-$th spot, we obtain vector fields $X_j$ defined on page \pageref{etot}.

Now we can define $d\pi_\lambda(\Xi)$ as a linear operator on $L^2(\R^n)$ determined by the rule
$$
d\pi_\lambda(\Xi)\f=\frac{d}{ds}\pi_\lambda(\gamma(s))\f\Big|_{s=0}\,.
$$

It is easy to see that 
$$
\aligned
d\pi_\lambda(X_{2k-1}) =\xi_k\\
d\pi_\lambda(X_{2k}) =\eta_k\\
d\pi_\lambda(Z) =\zeta\,,
\endaligned
$$
where, for $v\in\R^n$, $\xi_kf(v)=i\lambda v_kf(v)$, $\eta_kf(v)=\frac{\pd f}{\pd v_k}(v)$ and $\zeta f(v)=\frac i4\lambda f(v)$. 

It is possible to extend $d\pi_\lambda$ to 
the {\it universal enveloping algebra} of $\h^n$ (for definitions see \cite{Ta}, \cite{K} or \cite{Hu}), of which we can think as the unital associative algebra consisting of all left-invariant differential operators on $\cH^n$ with the binary operation being just the composition of operators. 
Consequently we get 
$$
d\pi_\lambda(\LL)=-\Delta+\lambda^2|v|^2\,,
$$
i.e. 
$d\pi_\lambda(\LL)$ is 
the (scaled) Hermite operator.  

By applying the method of transference \cite[Theorem 2.4]{CW} to the representation $\pi_\lambda$ and a single operator $X_{2k-1}\LL^{-1/2}$, one can  verify that we get 
$d\pi_\lambda(X_{2k-1}\LL^{-1/2})$, which is the same as $\xi_k L^{-1/2}$. Similarly, $X_{2k}\LL^{-1/2}$ is transferred to $\eta_k L^{-1/2}$.  
Hence if we transfer the operator $(X_{1}\LL^{-1/2},\hdots,X_{2n}\LL^{-1/2})$ the result is 
$(\xi_1 L^{-1/2},\eta_1 L^{-1/2},\hdots,\xi_n L^{-1/2},\eta_n L^{-1/2})$.
This is admittedly not the same as $2(\z_1 L^{-1/2},\z_1^* L^{-1/2},\hdots,\z_n L^{-1/2},\z_n^* L^{-1/2})$, but by \eqref{debeluh} the euclidean norms of these two vectors coincide when $\lambda=1$. Since $\pi_\lambda$ maps elements of $\He^n$ into contractions, we obtain, again by Theorem 2.4 from \cite{CW}, that $\nor{{\bf R}}_{B(L^p(\R^n))}\leqslant \nor{{\bf R}_{\He^n}}_{B(L^p(\He^n))}$. In words, the $L^p$ norm of the vector Hermite-Riesz transform is majorized by the $L^p$ norm of the vector Heisenberg Riesz transform. From \eqref{msu_notredame} we finally get the following: \vskip 7pt

{\it For every $\e>0$ there is $C_\e>0$ such that}
\begin{equation}
\label{predavanjeseminarka}
\nor{{\bf R}f}_p\leqslant C_\e p^{1+\e}\nor{f}_p  \,.
\end{equation}

\subsubsection*{Discussion: $1+\e$ or $1$}

We saw three results where the operators involved were estimated in $L^p$ by $O((p-1)^{-1-\e})$,  $p\rightarrow 1$, or $O(p^{1+\e})$, $p\rightarrow \infty$. These are: 1) estimates of the parabola Riesz transform $\mathcal{H}$, 2) the  estimate  \eqref{msu_notredame}
of $\nor{{\bf R}_{\He^n}f}_p$, and 3) the last estimate \eqref{predavanjeseminarka} of $\nor{{\bf R}f}_p$. 

The logic is that the estimate of $\mathcal{H}$ implies that of $\nor{{\bf R}_{\He^n}f}_p$, and this, in its turn, implies the estimate of $\nor{{\bf R}f}_p$.
On the other hand, we presented here Corollary \ref{RieszH}, where the estimate of $\nor{{\bf R}f}_p$ is obtained without $\e$, i.e. it is linear.   Two natural questions arise:

Q. 1: Is it possible to obtain a linear estimate of $\nor{{\bf R}_{\He^n}f}_p$?

Q. 2: Is it possible to obtain a linear estimate of $\|\mathcal{H}\|_p$?

We strongly believe that the answer to Q.1 is ``yes" and that the answer to Q.2 is ``no".

Let us comment on that. In what concerns Q.1 we believe that the Bellman function approach
used in the present paper is capable to treat a very wide range of Riesz transforms -- always giving estimates which are linear in $p$ and dimension-free.

As for
Q.2, it has been proved in \cite{TW} that a certain class of operators map
$L\log\log L$ into $L^{1,\infty}$. Parabola Riesz transform $\mathcal{H}$  belongs to this class, and now there is a strong feeling that this is sharp. In particular, $L\log^{\e}L$ is mapped to $L^{1,\infty}$, and, hence, $L\log^{1+\e}L$ to $L^1$, see \cite{BL}, \cite{TW}. But one cannot get rid of $\e$, because $L\log L$ is probably not mapped to $L^1$!
However, Yano's extrapolation theorem \cite{Y} 
implies that the linear estimate of $\mathcal{H}$ would give $L\log L$ to $L^1$ action. So we come to a ``contradiction".

\smallskip
Coming back to Q.1, it seems natural to try to prove an analogue of our bilinear embedding (Theorem \ref{bilinher}) in the context of
$\He^n$.  Then one can hope to have the analog of Corollary \ref{RieszH}.

\subsection{Sharpness of the linear estimate}

We believe that the linear $p-1$ estimate of Theorem \ref{bilinher} and its Corollary \ref{RieszH}  cannot be improved. There are several examples when similar singular operators got the linear estimate from below.
Each time it is a separate and non-trivial task. We believe this should be feasible and will be a subject of future efforts.

\bigskip

\noindent\hskip -15pt{\sc Oliver Dragi\v{c}evi\'c}, Faculty of Mathematics and Physics, University of Ljubljana, and Institute of Mathematics, Physics and Mechanics,\\ Jadranska 19, SI-1000 Ljubljana, Slovenia \\ {\tt oliver.dragicevic@fmf.uni-lj.si}\\

\noindent\hskip -15pt{\sc Alexander Volberg}, Department of Mathematics, Michigan State University, East Lansing, MI 48824, USA, and School of Mathematics, University of Edinburgh, Edinburgh EH9 3JZ, UK \\ {\tt volberg@math.msu.edu, A.Volberg@ed.ac.uk}


\end{document}